\input amstex
\documentstyle{amsppt}

\define\Cee{{\Bbb C}}

\define\Pee{\Bbb P}

\define\Zee{\Bbb Z}

\define\hil#1{{X^{[#1]}}}

\define\W{\widetilde}
\define\w{\tilde}

\def\exact#1#2#3{0\rightarrow{#1}\rightarrow{#2}\rightarrow{#3}\rightarrow0}

\define\im{\operatorname{Im}}
\define\Hom{\operatorname{Hom}}

\define\Sym{\operatorname{Sym}}

\define\Id{\operatorname{Id}}
\define\Pic{\operatorname{Pic}}

\define\Supp{\operatorname{Supp}}

\define\Hilb{\operatorname{Hilb}}

\define\Grass{\operatorname{Grass}}

\define\proof{\demo{Proof}}
\def\stop {\nobreak$\quad$\lower 1pt\vbox{
    \hrule
    \hbox to 7pt{\vrule height 7pt\hfil\vrule height 7pt}
      \hrule}\ifmmode\relax\else\par\medbreak \fi}
\define\endproof{\stop\enddemo}
\define\endstatement{\rm }
\define\theorem#1{\medskip\noindent {\bf {Theorem #1}.} \it}
\define\lemma#1{\medskip\noindent {\bf {Lemma #1}.} \it}
\define\proposition#1{\medskip\noindent {\bf {Proposition #1}.} \it}
\define\corollary#1{\medskip\noindent {\bf {Corollary #1}.} \it}
\define\claim#1{\medskip\noindent {\bf {Claim #1}.} \it}

\define\section#1{\bigskip\noindent {\bf #1}}
\define\ssection#1{\medskip\noindent{\bf #1}}

\def\mapdown#1{\Big\downarrow\rlap{$\vcenter{\hbox{$\scriptstyle#1$}}$}}
\def\mapright#1{\smash{   \mathop{\longrightarrow }\limits^{#1}}}

\def\mapdown#1{\Big\downarrow
    \rlap{$\vcenter{\hbox{$\scriptstyle#1$}}$}}

\def\today{\ifcase \month \or January\or February\or  
  March\or April\or May\or June\or July\or August\or 
  September\or October\or November\or December\fi\   
  \number \day, \number \year}

\loadbold
\hoffset=0.3 in
\voffset=0.6 in

\leftheadtext{Wei-ping Li, Zhenbo Qin and  Qi Zhang}
\rightheadtext{Geometry of Hilbert schemes of points}

\topmatter
\title
On the geometry of the Hilbert schemes of points in the projective
plane
\endtitle
\author {Wei-ping Li$^1$, Zhenbo Qin$^2$ and Qi Zhang$^3$}
\endauthor

\thanks
${}^1$Partially supported by the grant HKUST6170/99P
\endthanks
\thanks 
${}^2$Partially supported by an NSF grant and 
an Alfred P. Sloan Research Fellowship
\endthanks
\thanks 
${}^3$Partially supported by an NSA grant and 
MU Research Board
\endthanks
\abstract
We determine the nef cone and effective cone of the Hilbert scheme $\hil{n}$
of points on the projective plane $X = \Pee^2$,
characterize all the rational curves in $\hil{n}$ of degree-$1$
with respect to some very ample line bundle, 
describe the moduli spaces of these curves, 
and study the contraction of the extremal ray on $\hil{n}$.
\endabstract

\endtopmatter
\noindent ${}^{ \,\,\,\,\, 1}$Department of Mathematics, HKUST, Clear Water Bay,
Kowloon, Hong Kong

\noindent ${}^{2,3}$Department of Mathematics, University of Missouri, Columbia,
MO 65211, USA

\NoBlackBoxes
\TagsOnRight
\document

\section{1. Introduction}

The Hilbert scheme $\hil{n}$ of points in
a smooth projective surface $X$
is a desingularization of the $n$-th symmetric product 
$\Sym^n(X)$ of $X$ (see \cite{Fo1}). An element $\xi$ in $\hil{n}$ 
is a length-$n$ $0$-dimensional closed subscheme of $X$. 
Recently, there are two surprising discoveries, mainly due to the work 
of G\"ottsche \cite{Got}, Nakajima \cite{Na1} and Grojnowski \cite{Gro}, 
that the Hilbert schemes $\hil{n}$ ($n \ge 0$) have relationships
with modular forms on the one hand and with representations of
infinite dimensional Lie algebras on the other hand.
Their work basically deals with the cohomology groups
(the Hodge numbers as well) of $\hil{n}$.
The geometry of $\hil{n}$ is another interesting
area of research. An excellent example in this regard is 
Beauville's work \cite{Bea} which says that $\hil{n}$ admits 
a hyperk\"ahler structure when $X$ is a $K3$ surface. In this paper, 
we shall study the geometry of the Hilbert scheme $\hil{n}$ 
when $X$ is the projective plane $\Bbb P^2$. To be more specific, 
we shall determine the nef cone and effective cone of $\hil{n}$,
characterize all the rational curves in $\hil{n}$ of degree-$1$
with respect to some very ample line bundle on $\hil{n}$, 
and study the contraction of the extremal ray on  $\hil{n}$.
Throughout the paper, we assume that $n \ge 2$.

Let $X=\Bbb P^2$, $\ell$ be a line in $X$, 
and $x_1,\ldots, x_{n-1} \in X$ be distinct but 
fixed points with $x_i \notin \ell$ for $i=1, \ldots, (n-1)$. Let 
$M_2(x_1)=\{\xi\in \hil{2}\,|\, \Supp(\xi)=\{x_1\} \}$
be the punctual Hilbert scheme parameterizing length-$2$ 
$0$-dimensional subschemes supported at $x_1$. 
It is known that $M_2(x_1) \cong \Bbb P^1$. 
We define two curves in $\hil{n}$:
$$\align
 \beta_{\ell}
&=\big\{\,x+x_1+\ldots+x_{n-1} \in \hil{n}\,|\,x\in \ell \big\},\tag1.1\\
 \beta_n
&=\big\{\, \xi+ x_2+ \ldots+x_{n-1} \in \hil{n}\,|\,\xi\in M_2(x_1) \big\}\tag 1.2
\endalign$$
Both $\beta_n$ and $\beta_{\ell}$ are smooth rational curves on $\hil{n}$. 
It can be shown that the $1$-cycle $\beta_{\ell}-(n-1)\beta_n$ 
in the Chow group $A_1(\hil{n})$ is rationally equivalent to 
a smooth rational curve. 
In section 3, we determine the effective cone of $\hil{n}$.

\theorem{A} The effective cone of $\hil{n}$ is spanned by 
$\beta_n$ and $\beta_{\ell} -(n-1)\beta_n$.
\endstatement

Next we introduce the following two divisors on 
the Hilbert scheme $\hil{n}$:
$$
B_n     =\big\{ \, \xi\in \hil{n}\,|\, |\Supp(\xi)|<n\, \big\}
\quad\hbox{and}\quad
D_{\ell}=\big\{ \, \xi \in \hil{n}\, |\, \Supp(\xi)\cap \ell 
         \neq\emptyset\big\}.
$$
It can be shown that $(nD_{\ell}-B_n/2)$ is very ample, and
$\beta_n \cdot (nD_{\ell}-B_n/2)=1 =
(\beta_{\ell}-(n-1)\beta_n)\cdot (nD_{\ell}-B_n/2)$.
Moreover, curves of degree-$1$ with respect to $(nD_{\ell}-B_n/2)$ 
are homologous to either $\beta_n$ or $\beta_{\ell}-(n-1)\beta_n$.
We are interested in determining all these minimal degree curves.

In section 4, we characterize all the curves $\gamma$  
homologous to $\beta_n$. To state the result, we put  
$\hil{n, n+1}=\big\{ \, (\xi, \eta)\in \hil{n}\times \hil{n+1}\,
|\, \xi\subset \eta \,\big\}$, and define
$${f_{n+1}}\colon \hil{n, n+1}\rightarrow
\hil{n},\qquad
{\psi_{n+1}}\colon\hil{n, n+1}\rightarrow\Cal Z_{n+1}
\subset \hil{n+1}\times X$$
where $f_{n+1}$ is the restriction of the projection  
$\hil{n}\times \hil{n+1} \to \hil{n}$, $\Cal Z_{n+1}$ is 
the universal  subscheme of $\hil{n+1}\times X$, 
and $\psi_{n+1}(\xi, \eta) =
\big(\eta, \Supp(I_\xi/I_\eta) \big)\in \Cal Z_{n+1}$.
It is known (see \cite{ES2}) that $\psi_{n+1}: \hil{n, n+1} \to \Cal Z_{n+1}$ 
is canonically isomorphic to
$\Pee(\omega_{\Cal Z_{n+1}})$ where $\omega_{\Cal Z_{n+1}}$ is 
the dualizing sheaf of $\Cal Z_{n+1}$. 
So the fiber $(\psi_{n+1})^{-1}(\eta, x)$ over 
$(\eta, x) \in \Cal Z_{n+1}$ is isomorphic to the projective space
$\Pee(\omega_{\eta_x} \otimes \Cee_x)$. 

\theorem{B} A curve $\gamma$ in $\hil{n}$ is homologous to 
$\beta_n$ if and only if $\gamma=f_{n+1}(C)$ where $C$ is a line
in the projective space $(\psi_{n+1})^{-1}(\eta, x)$ 
with $(\eta, x) \in \Cal Z_{n+1}$.
Moreover, in this case, the point $(\eta, x)$ and the line $C$
are uniquely determined by $\gamma$.
\endstatement

We shall show that the moduli space of all these curves 
has dimension $(2n-2)$ and its top stratum consists of all the curves
$\gamma$ of the form (1.2). Moreover, the normal bundles of 
these curves $\gamma$ in $\hil{n}$ are determined.

In section 5, we characterize all the curves $\gamma$ homologous to 
$\beta_{\ell} -(n-1)\beta_n$. Recall that $\Hilb^n(\Pee^1) \cong 
\Sym^n(\Pee^1) \cong \Pee^n$. So $\Hilb^n(C) \cong \Pee^n$
for every line $C$ in $X=\Bbb P^2$.

\theorem{C} A curve $\gamma$ in $\hil{n}$ is homologous to 
$\beta_{\ell}-(n-1)\beta_n$ if and only if there exists a line $C$ in 
$X=\Bbb P^2$ such that $\gamma$ is a line in $\Hilb^n(C)\subset \hil{n}$.
Moreover, in this case, the line $C$ in $X=\Bbb P^2$ is uniquely
determined by the curve $\gamma$.
\endstatement

Note that the set of all lines in $X=\Pee^2$ is the dual space $(\Pee^2)^*$,
and $\Pee((T_{(\Pee^2)^*})^*)$ is isomorphic to the universal family
$\{ (x, C) \in \Pee^2 \times (\Pee^2)^*| \, x \in C \}$.
We prove that the set of all the curves $\gamma$
in $\hil{n}$ homologous to $\beta_{\ell}-(n-1)\beta_n$ is
parameterized by the Grassmannian bundle
$\Grass(\Sym^n((T_{(\Pee^2)^*})^*),2)$ over $(\Pee^2)^*$
(here we refer to Conventions below for 
our convention on the Grassmannians).

In section 6, we study a morphism $\varphi_1$ on $\hil{n}$ which is 
defined as follows. 
Let $\text{Grass}\big(H^0\big(X, \Cal O_X(n-1)\big),n\big)$
be the Grassmannian of $n$-dimensional quotients of the vector space
$H^0\big(X, \Cal O_X(n-1)\big)$. By a result of 
Beltrametti and Sommese \cite{BSG} on $k$-very ampleness, 
there exists a morphism 
$\varphi_1\colon\, \hil{n} \to
\text{Grass}\big(H^0\big(X, \Cal O_X(n-1)\big), n\big)$
mapping $\xi \in \hil{n}$ to the $n$-dimensional quotient space 
$H^0\big(X, \Cal O_{\xi}\otimes \Cal O_X(n-1)\big)$ of 
$H^0\big(X, \Cal O_X(n-1)\big)$.
We show that 
when $n \ge 3$, $\varphi_1$ is the contraction of an extremal ray.
Moreover, $\varphi_1: \hil{n} \to \varphi_1(\hil{n})$ 
is a small contraction when $n \ge 4$.

\theorem{D} When $n \ge 4$, the flip for the small contraction
$\varphi_1$ exists.
\endstatement

\medskip
All the results (including Theorem B) in section 4 are proved 
for a simply-connected smooth projective surface $X$.
They have been used to compute 
the $1$-point Gromov-Witten invariants of the Hilbert scheme $\hil{n}$, 
and to study the quantum cohomology of $\hil{n}$. 
We refer to \cite{L-Q} for details. 

\medskip\noindent
{\bf Conventions:} 
Let $0 \le k \le n$ and $V$ be an $n$-dimensional vector space. 
We use the Grassmannian $\Grass(V, k)$ to denote the set of 
all $k$-dimensional quotients of $V$, or equivalently, 
the set of all $(n-k)$-dimensional subspaces of $V$.
Also, we take $\Bbb P(V) = \Grass(V, 1)$. 
So the set of lines in $\Bbb P(V)$ is the Grassmannian $\Grass(V, 2)$.

\medskip\noindent
{\bf Acknowledgments:} 
The authors thank D. Edidin, L. Ein,  R. Friedman, S. Katz, 
Y. Kawamata and Y. Ruan for 
stimulating discussions. 
In addition, the first author thanks S. Katz
and the Mathematics Department of Oklahoma State University for their
hospitality during his visit there in June of 1999.
This work is partially supported by a grant from
the American Institute of Mathematics. 

\section{2. Basic facts about the Hilbert scheme of points in a surface}

In this section, we shall fix some notations, 
and review some basic facts for the Hilbert scheme
of points in a smooth projective surface.

Let $X$ be a simply-connected smooth projective surface,
and $\hil{n}$ be the Hilbert scheme of points in $X$.
An element in $\hil{n}$ is represented 
by a length-$n$ $0$-dimensional closed subscheme $\xi$ of $X$. 
For $\xi \in \hil{n}$, let $I_{\xi}$ and 
$\Cal O_\xi$ be the corresponding sheaf of ideals and 
structure sheaf respectively.
For a point $x \in X$, let $\xi_x$ be the component of $\xi$ 
supported at $x$ and $I_{\xi, x} \subset \Cal O_{X, x}$ be 
the stalk of $I_{\xi}$ at $x$. 
It is known from \cite{Fo1} that $\hil{n}$ is smooth. 
In $\hil{n}\times X$, we have the universal codimension-$2$ subscheme:
$$\Cal Z_n=\{(\xi, x) \subset \hil{n}\times X \, |\,
x\in \Supp{(\xi)}\}\subset \hil{n}\times X.\tag 2.1$$
Define the incidence variety
$\hil{n-1, n}=\{ (\xi,\eta)\in \hil{n-1}\times \hil{n} \,|\,
\xi \subset \eta \} \subset \hil{n-1}\times \hil{n}$.
It is known that $\hil{n-1,n}$ is smooth and of dimension $2n$.
Let $X^n$ and $X^{(n)}=\text{Sym}^n(X)$ be the $n$-th Cartesian product
and symmetric product respectively.
Then, we have the Hilbert-Chow morphism $\rho_n$ and
the quotient map $\upsilon_n$:
$$\rho_n\colon\,\, \hil{n} \to X^{(n)}, \qquad  
\upsilon_n\colon\,\, X^n \to  X^{(n)}.    \tag2.2$$

For a subset $Y \subset X$, we define the subset $M_n(Y)$ in 
the Hilbert scheme $\hil{n}$:
$$M_n(Y) = \{ \xi \in \hil{n}| \Supp(\xi) \text{ is a point in } Y \}
\subset \hil{n}.  \tag 2.3$$
So for $x \in X$,
$M_n(x)$ is just the punctual Hilbert scheme of points on $X$ at $x$.
It is known that  $M_n(x)$ are 
isomorphic for all  $X$ and all  $x \in X$.
Let $\Hilb^n(\Bbb C^2, 0)$ be 
the punctual Hilbert scheme of points on $\Bbb C^2$ at the origin.
We adopt the isomorphism 
$$M_n(x) \cong \Hilb^n(\Bbb C^2, 0).\tag 2.4$$

Given $Y \subset \hil{n-k}$ and $\eta \in \hil{k}$
with $\{\Supp(\xi)|{\xi\in Y}\}\cap \Supp(\eta) =\emptyset$, 
we use $Y+\eta$ to represent the subset in $\hil{n}$ consisting 
of all the points $\xi + \eta$ with $\xi \in Y$.

\medskip\noindent
{\bf Notation 2.5.} (see \cite{ES2}) Let $f_n$ and $g_n$
be the projections of $\hil{n-1, n} \subset \hil{n-1}\times \hil{n}$ 
to $\hil{n-1}$ and $\hil{n}$ respectively. 
Let $p_n$ and $q_n$ be the projections of 
$\hil{n}\times X$ to $\hil{n}$ and $X$ respectively.
Define the morphisms $\pi_n\colon \,\, \Cal Z_n\to \hil{n}$ by 
$\pi_n(\eta, x)=\eta$ and $q\colon \,\, \hil{n-1, n}\to X$ by 
$q(\xi, \eta)= \Supp(I_{\xi}/I_{\eta})$. In addition, 
we put $\psi_n =(g_n, q) \colon \,\, \hil{n-1, n}\to \Cal Z_n$
and $\phi_n=(f_n, q)\colon \,\, \hil{n-1, n}\to \hil{n-1}\times X$.

Next, we review some results on the homology groups of 
$\hil{n}$ due to G\"ottsche \cite{Got}, Grojnowski \cite{Gro},
and Nakajima \cite{Na1}. 
Let $H=H_0(X,\Bbb C)\oplus H_2(X, \Bbb C)\oplus H_4(X,\Bbb C)$, 
and $S=S(H\otimes_{\Bbb C}t^{-1}\Bbb C[t^{-1}])$ be
the symmetric algebra over the infinite dimensional complex vector space
$H\otimes_{\Bbb C}t^{-1}\Bbb C[t^{-1}]$. 
For $\alpha \in H_i(X, \Bbb C)$ and $k > 0$, let $|\alpha| = i$ 
and associate to $\alpha \otimes t^{-k}$ a bidegree $(k, 2k-2+|\alpha|)$.
The results in \cite{Got, Gro, Na1} state that 
$\displaystyle{S \cong \bigoplus\limits^{+ \infty}_{n=0}
\bigoplus\limits^{4n}_{k=0} H_k(\hil{n}, \Bbb C)}$.
This identification can be understood as follows.
For $i = 1, \ldots, k$, let $\alpha_i \in H_{|\alpha_i|}(X, \Bbb C)$ 
be represented by a real submanifold $A_i$ such that
$A_1, \ldots, A_k$ are in general position.
Then, the element $(\alpha_1\otimes t^{-n_1})
\ldots(\alpha_k\otimes t^{-n_k}) \in S$ has bidegree
$$\left (\sum_{i=1}^k n_i, \sum\limits_{i=1}^k(2n_i-2+|\alpha_i|)
\right ).    \tag 2.6$$
Moreover, up to a scalar multiple, $(\alpha_1\otimes t^{-n_1})
\ldots(\alpha_k\otimes t^{-n_k})$
is represented by the closure of
the real-$\sum_{i=1}^k(2n_i-2+|\alpha_i|)$-dimensional subset:
$$\{ \xi_1+\ldots + \xi_k \in \hil{n}|
\xi_i \in M_{n_i}(A_i) \text{ and } \Supp(\xi_i) \cap \Supp(\xi_j) =
\emptyset \text{ for } i \ne j \}  \tag 2.7$$
where $n=\sum_{i=1}^k n_i$ and $M_{n_i}(A_i)$ is the subset of
$\hil{n_i}$ defined by (2.3).

In particular, if $\alpha \in H_2(X, \Bbb C)$ is represented by
a real surface $C$ in $X$ and if 
$D_C {\overset \text{def} \to =}
\{\xi\in \hil{n} \,|\, C\cap \Supp(\xi)\neq \emptyset\},$
then we obtain the homologous relation:
$$D_C\sim {1/(n-1)!} \cdot  (\alpha \otimes t^{-1})
(X\otimes t^{-1})^{n-1} \in S.  \tag 2.8$$
Note that if $C$ is an algebraic curve in $X$,
then $D_C$ is a divisor in $\hil{n}$. In fact, 
$\Cal O_{\hil{n}}(D_C) = (\rho_n)^*\Cal L_C$ 
where $\Cal L_C$ is the unique line bundle over $X^{(n)}$ satisfying
$(\upsilon_n)^*\Cal L_C = \otimes_{i=1}^n (\lambda_i)^*\Cal O_X(C).$
Here $\lambda_i$ is the projection of $X^n$ to the $i$-th factor. 
Thus if $n\ge 2$ and the curve $C$ is an ample divisor on $X$, 
then $D_C$ is nef and big, but not ample.
By linearity, we extend the notation $D_C$ to any divisor $C$ on $X$.

When $n \ge 2$, another distinguished divisor on $\hil{n}$ is defined by
$$B_n {\overset \text{def} \to =} \{\xi\in \hil{n}\,|\, |\Supp(\xi)|<n\} 
\sim {1/(n-2)!} \cdot 
(X\otimes t^{-2})(X\otimes t^{-1})^{n-2}. \tag 2.9$$

\lemma{2.10} Let $n \ge 2$ and  $X$ be simply-connected.
Let real surfaces $\{C_1, \ldots, C_s \}$ form  a basis of
$H_2(X,\Bbb C)$. Then $B_n, D_{C_1}, \ldots, D_{C_s}$ form 
a basis of $H_{4n-2}(\hil{n}, \Bbb C)$.
\endstatement
\proof 
Let $\alpha_i$ be the homology class represented by $C_i$.
Expand the basis $\{\alpha_1, \ldots, \alpha_s \}$ of
$H_2(X, \Bbb C)$ to the basis
$\{\alpha_0 = x, \alpha_1, \ldots, \alpha_s,\alpha_{s+1}=X \}$
of $H = H_0(X, \Bbb C) \oplus H_2(X, \Bbb C) \oplus H_4(X, \Bbb C)$.
By (2.6), a basis of $H_{4n-2}(\hil{n}, \Bbb C)$
consists of elements of the form
$(\alpha_{m_1} \otimes t^{-n_1})\ldots(\alpha_{m_k}\otimes t^{-n_k})$
with $n_i\ge 1$, $\sum_{i=1}^k n_i=n$,
and $\sum_{i=1}^k(2n_i-2+|\alpha_{m_i}|)=4n-2$.
Hence, $\sum_{i=1}^k(|\alpha_{m_i}|-2)=2n-2$.
Since $|\alpha_{m_i}| \in \{0,2,4\}$, 
either $|\alpha_{m_i}|=4$ for all $i$ 
or $|\alpha_{m_{i_0}}|=2$ for some $i_0$
and $|\alpha_{m_j}|=4$ for $j \ne i_0$.

In the first case, $\alpha_{m_i}=X$ for all $i$.
Since $2n-2=\sum_{i=1}^k(4-2)$, $k=n-1$.
Since $n=\sum_{i=1}^k n_i$,
$n_1=2$ and $n_2 = \ldots = n_{n-1}=1$
up to an ordering of $n_1, \ldots, n_{n-1}$. So 
$(\alpha_{m_1}\otimes t^{-n_1})\ldots(\alpha_{m_k}\otimes t^{-n_k})
\sim (n-2)! \cdot B_n$ by (2.9).
In the second case, we may let $i_0=1$. 
Then $2n-2=\sum_{i=1}^k(|\alpha_{m_i}|-2)=2(k-1)\le 2n-2$.
So $k=n$, $n_i=1$ for all $i$, and 
$(\alpha_{m_1} \otimes t^{-n_1})\ldots(\alpha_{m_k} \otimes t^{-n_k})  
\sim (n-1)! \cdot D_{C_{m_1}}$ by (2.8).
\endproof

\section{3. The nef cone and effective cone of $X^{[n]}$ when $X=\Bbb P^2$}

In this section, we review some construction in \cite{BSG}
and use it to determine the nef cone and effective cone of $\hil{n}$
when $X$ is the projective plane $\Bbb P^2$.

\medskip\noindent
{\bf Definition 3.1.} (see \cite{BSG}) 
A line bundle $L$ on a projective surface $X$ is 
{\it $k$-very ample} if the restriction 
$H^0(X, L) \to H^0(X, \Cal O_{\xi}\otimes L)$
is surjective for every $\xi\in \hil{k+1}$.
\medskip

The concept of $k$-very ampleness relates  
$\hil{k+1}$ to a Grassmannian as follows.
The surjective map in Definition 3.1 represents
an element in $\Grass(H^0(X, L), k+1)$.
So we get a morphism
$\varphi_k(L) \colon\,\, \hil{k+1} \to \Grass(H^0(X, L), k+1)$
when $L$ is $k$-very ample. Letting $h = h^0(X, L)$ and 
$\frak P \colon\,\, \Grass(\Cee^h, k+1)
\to \Bbb P((\wedge^{h-(k+1)}\Bbb C^h)^*)$ be the Pl\"ucker embedding,
then we see from the Appendix of \cite{BSG} that 
$$(\frak P \circ \varphi_k(L))^* \Cal H
=\Cal O_{\hil{k+1}}\big(D_{c_1(L)}-B_{k+1}/2 \big)  \tag 3.2$$
where $\Cal H$ is the hyperplane line bundle over the projective space
$\Bbb P((\wedge^{h-(k+1)} \Bbb C^h)^*)$.

\smallskip\noindent
{\bf Remark 3.3.} Let $k\ge 1$ be an integer. 
By the Theorem 2.1 in \cite{BSG},
there exists a very ample divisor $C$ on $X$
such that the line bundle $L=\Cal O_X(C)$ is $(k+1)$-very ample.
Then by the Main Theorem in \cite{C-G}, the morphism
$\varphi_{k}(L)$ is an embedding. By (3.2), 
the divisor $(D_C-B_{k+1}/2)$ on $\hil{k+1}$ is very ample. 
\medskip 

In the rest of this section, we assume that $X=\Bbb P^2$.

\lemma{3.4} Let $X=\Bbb P^2$. Then the line bundle $\Cal O_X(n)$
is $n$-very ample.
\endstatement
\proof We use the Theorem 2.1 in \cite{BSG}.
Let $\ell$ be a line on $\Bbb P^2$, $C=(n+3)\ell$, and $k=n$.
Then, $K_X+C=n\ell$ and $C\cdot C=(n+3)^2 > 4n+5$.
Suppose $\Cal O_X(n)=\Cal O_X(K_X+C)$ is not $n$-very ample.
By the Theorem 2.1 in \cite{BSG}, 
there exists an effective divisor $D$ on $X$ such that  
$C\cdot D-n-1\le D\cdot D < C\cdot D/2< n+1$.
Let $D \sim d\ell$. Then, $d \ge 0$ and
$(n+3)d-n-1\le d^2<(n+3)d/2< n+1$. Since $d^2<(n+3)d/2$, $d>0$.
Now, $(n+3)d/2<(n+1)$ implies $d<2$. So $d=1$
which contradicts to $(n+3)d-n-1\le d^2$.
Hence, $\Cal O_X(n)$ is $n$-very ample.
\endproof

Note that for $k \ge 1$, $k$-very ampleness implies $(k-1)$-very ampleness. 
Therefore, by Lemma 3.4, we get the following two morphisms:
$$\align
\varphi_1 {\overset \text{def} \to =} \varphi_{n-1}\big(\Cal O_X(n-1)\big)
 &\colon\,\, \hil{n} \to
  \Grass\big(H^0\big(X, \Cal O_X(n-1)\big), n\big),\tag 3.5 \\
\varphi_2 {\overset \text{def} \to =} \varphi_{n-1}\big(\Cal O_X(n)\big)
 &\colon\,\, \hil{n}\to
  \Grass\big(H^0\big(X, \Cal O_X(n)\big), n\big). \tag 3.6
\endalign$$ 
For $i=1$ or $2$, let $\ell_{i}=h^0\big(X, \Cal O_X(n-2+i)\big)\
={n+i\choose 2}$, $N_i = \dim (\wedge^{\ell_i-n}\Bbb C^{\ell_i})-1$,
$\frak P_i\colon\,\, \Grass(\Bbb C^{\ell_i}, n) \to \Bbb P^{N_i}$ 
be the Pl\"ucker embedding, $\Phi_i \,\,{\overset \text{\rm def} \to =}\,\,
\frak P_i \circ \varphi_i \colon\,\, \hil{n} \to \Bbb P^{N_i}$,
and $\Cal H_i$ be the hyperplane line bundle over $\Bbb P^{N_i}$. By (3.2),  
$$(\Phi_1)^*\Cal H_1
 =\Cal O_{\hil{n}}\big((n-1)D_{\ell}-B_n/2 \big),\quad
(\Phi_2)^*\Cal H_2
 =\Cal O_{\hil{n}}(nD_{\ell}-B_n/2).\tag 3.7$$

\lemma{3.8} Let $X = \Pee^2$. Then the morphism $\varphi_2$ in (3.6)
is an embedding of the Hilbert scheme $\hil{n}$,
and the divisor $(nD_{\ell}-B_n/2)$ is very ample.
\endstatement
\proof
By Lemma 3.4, $\Cal O_X(n)$ is $n$-very ample.
So by the Main Theorem in \cite{C-G}, $\varphi_2$ is an embedding.
By (3.7), $(nD_{\ell}-B_n/2)$ is very ample.
\endproof

Next, we study $((n-1)D_{\ell}-B_n/2)$.
We begin with a construction. Let $\xi = x_1+ \ldots+ x_{n-1} \in \hil{n-1}$ 
where $x_1, \ldots, x_{n-1}$ are distinct points in the line $\ell \subset X$. 
Let $(\ell+\xi)$ be the closure of 
$\big(\ell-\Supp(\xi)\big)+\xi$ in $\hil{n}$. Alternatively, consider:
$$\matrix 
    \w \ell_{\xi} &\subset &\W X_{\xi}&\subset&\hil{n-1, n}&\mapright{g_n}&\hil{n}\\
    \mapdown{}    &&\mapdown{}&&\mapdown{\phi_n}\\
    \xi\times \ell&\subset &\xi\times X&\subset& \hil{n-1}\times X
\endmatrix  \tag 3.9
$$
where $\w \ell_{\xi}$ and $\W X_{\xi}$ are the strict transforms of 
$\xi\times \ell$ and $\xi\times X$ in $\hil{n-1, n}$ respectively.
Since $\phi_n$ is the blowup of $\hil{n-1}\times X$ along $\Cal Z_{n-1}$ 
(see the Proposition 2.2 of \cite{ES2}),
$\W X_{\xi}$ is isomorphic to the blowup of $\xi\times X \cong X$ 
at $x_1,\ldots, x_{n-1}$.
Let $E_i$ be the exceptional divisor in $\W X_{\xi}$ over $x_i$.
Then, $(\phi_n|_{\W X_{\xi}})^*(\xi\times \ell) =\w \ell_{\xi}+
\sum_{i=1}^{n-1}E_i$ in $A_1(\W X_{\xi})$.
Notice that $g_n(\w \ell_{\xi})=(\ell+\xi)$ and 
$g_n(E_i)=M_2(x_i)+x_1\ldots+x_{i-1}+x_{i+1}+\ldots+x_{n-1}$.
In fact, since $g_n|_{\W X_{\xi}}: \W X_{\xi}
\to g_n({\W X_{\xi}})$ is an isomorphism, we have
$$(g_n|_{\W X_{\xi}})_*(\w \ell_{\xi}) = (\ell+\xi) \quad \text{and} \quad
(g_n|_{\W X_{\xi}})_*(E_i) =\beta_n \tag 3.10$$
in $A_1(\hil{n})$. Here $\beta_n$ is defined in (1.2).
Also, recall $\beta_\ell$ from (1.1).

\lemma{3.11} With the above notations, 
$(\ell+\xi) = \beta_\ell-(n-1)\beta_n$ in $A_1(\hil{n})$.
\endstatement
\proof
Choose another line $\ell'$ such that $\Supp(\xi) \cap \ell' = \emptyset$. 
Then in $A_1(\hil{n})$, 
$$\align 
(g_n|_{\W X_{\xi}})_*(\phi_n|_{\W X_{\xi}})^*(\xi\times \ell)
 &=(g_n|_{\W X_{\xi}})_*(\phi_n|_{\W X_{\xi}})^*(\xi \times \ell')
  =(\ell'+\xi)=\beta_{\ell'} = \beta_\ell,  \\
(g_n|_{\W X_{\xi}})_*(\phi_n|_{\W X_{\xi}})^*(\xi\times \ell)
 &=(g_n|_{\W X_{\xi}})_*(\w \ell_{\xi}+\sum_{i=1}^{n-1}E_i)
  =(\ell+\xi)+(n-1)\beta_n
\endalign$$
where we have used (3.10). 
So $(\ell+\xi)=\beta_\ell-(n-1)\beta_n$ in $A_1(\hil{n})$.
\endproof

\proposition{3.12} Let $X=\Pee^2$,
and $\varphi_1$ be the morphism defined in (3.5).

{\rm (i)} If $n\ge 2$, then the divisor $((n-1)D_{\ell}-B_n/2)$ 
on $\hil{n}$ is nef but not ample.

{\rm (ii)} If $n=2$, then $\varphi_1: \hil{2} \rightarrow
\Grass\big(H^0\big(X, \Cal O_X(1)\big), 2\big)=(\Bbb P^2)^*$
is isomorphic to the relative Hilbert scheme
$\Hilb^2(\Pee((T_{(\Pee^2)^*})^*)/(\Bbb P^2)^*)$ over 
the dual space $(\Bbb P^2)^*$.
\endstatement
\proof
(i) We have the following intersection numbers:
$$\beta_{\ell}\cdot D_{\ell}=1,\quad
\beta_{\ell}\cdot (-B_n/2)=0,\quad
\beta_n \cdot D_{\ell}=0,\quad
\beta_n \cdot (-B_n/2)=1 \tag 3.13$$
where  the last one follows from the Claim 2 in the Appendix of \cite{BSG}.
It follows that $\big((n-1)D_{\ell}-B_n/2\big)\cdot 
(\beta_{\ell}-(n-1)\beta_n) =0$.
By Lemma 3.11, $(\beta_{\ell}-(n-1)\beta_n)$ is rationally equivalent to
an effective curve. So $((n-1)D_{\ell}-B_n/2)$ is not ample.
By (3.7), the divisor $((n-1)D_{\ell}-B_n/2)$ is nef. 

(ii) When $n=2$, $\varphi_1=\varphi_1\big(\Cal O_X(1)\big)$: 
$\hil{2}\to \Grass\big(H^0\big(X, \Cal O_X(1)\big), 2\big)
=(\Bbb P^2)^*$.
An element $\xi\in \hil{2}$ determines a unique line
$\ell_{\xi}$ in $X$ and $\varphi_1(\xi)=\ell_{\xi}\in (\Bbb P^2)^*$.
So $\varphi_1$ is a surjection
whose fiber over a line $\ell \in (\Bbb P^2)^*$ is
$\Hilb^2(\ell) \cong \Sym^2(\ell) \cong \Bbb P^2$.
Note that 
$\{(x, \ell) \in \Bbb P^2 \times (\Bbb P^2)^*|\, \, x \in \ell \}
\cong \Pee((T_{(\Pee^2)^*})^*)  $.
So $\varphi_1: \hil{2} \to (\Bbb P^2)^*$ is isomorphic to
the relative Hilbert scheme 
$\Hilb^2(\Pee((T_{(\Pee^2)^*})^*)/(\Bbb P^2)^*)$ over $(\Bbb P^2)^*$.
\endproof

\theorem{3.14} Let $n \ge 2$ and $\ell$ be a line in $X = \Pee^2$. 

{\rm (i)} The nef cone of $\hil{n}$ is the cone spanned by
$D_{\ell}$ and $((n-1)D_{\ell}-B_n/2)$;

{\rm (ii)} The effective cone of $\hil{n}$ is the cone spanned by
$\beta_n$ and $(\beta_{\ell}-(n-1)\beta_n)$;

{\rm (iii)} Let $\beta$ be the homology class
of an effective curve in the Hilbert scheme $\hil{n}$. Then,
$\beta \cdot (nD_{\ell}-B_n/2) = 1$ if and only if  
either $\beta = \beta_n$ or $\beta = (\beta_{\ell}-(n-1)\beta_n)$.
\endstatement
\proof
(i) Note that $D_{\ell}$ is nef and big, but not ample
(see the paragraph containing (2.8)). By Proposition 3.12,
$((n-1)D_{\ell}-B_n/2)$ is also nef but not ample.
By the Theorem 1.1 (i) in \cite{ES1} and the results in \cite{Fo2, BSG}, 
the divisors $(-B_n/2)$ and $D_{\ell}$ form a basis for 
$A_{2n-1}(\hil{n}) \cong H_{4n-2}(\hil{n}, \Bbb Z)$.
So the nef cone of the Hilbert scheme $\hil{n}$ is the cone spanned by 
$D_{\ell}$ and $((n-1)D_{\ell}-B_n/2)$.

(ii) Since $(-B_n/2)$ and $D_{\ell}$ form a basis for 
$H_{4n-2}(\hil{n}, \Bbb Z)$, we see from (3.13) and the Poincar\'e duality
that $\{ \beta_{\ell}, \beta_n \}$ is a basis for 
$H_2(\hil{n}, \Zee)$. By (1.2) and Lemma 3.11, both $\beta_n$ and 
$(\beta_{\ell}-(n-1)\beta_n)$ are numerical equivalence classes of 
effective curves. By (3.13) again, $\beta_n \cdot D_\ell = 0, 
(\beta_{\ell}-(n-1)\beta_n) \cdot D_\ell =1$,
$\beta_n \cdot ((n-1)D_{\ell}-B_n/2) = 1, 
(\beta_{\ell}-(n-1)\beta_n) \cdot ((n-1)D_{\ell}-B_n/2) = 0$.
It follows from (i) that the effective cone of $\hil{n}$
is spanned by $\beta_n$ and $(\beta_{\ell}-(n-1)\beta_n)$.

(iii) Let $\beta = a\beta_n + b(\beta_{\ell}-(n-1)\beta_n)$.
Intersecting $\beta$ with the nef divisors $D_\ell$ and 
$((n-1)D_{\ell}-B_n/2)$, we see that $a$ and $b$ are nonnegative integers.
Since $(nD_{\ell}-B_n/2)$ is very ample and 
$\beta_n \cdot (nD_{\ell}-B_n/2) 
= (\beta_{\ell}-(n-1)\beta_n) \cdot (nD_{\ell}-B_n/2)=1$,
we conclude that $\beta \cdot (nD_{\ell}-B_n/2) = 1$ if and only if 
$\beta = \beta_n$ or $\beta = (\beta_{\ell}-(n-1)\beta_n)$.
\qed

\section{4. Curves homologous to $\beta_n$}

In this section, $X$ stands for an arbitrary simply-connected
smooth projective surface (not just $\Pee^2$).
We shall determine all the curves in the Hilbert scheme $\hil{n}$
homologous to $\beta_n$, and study their moduli space.

First of all, we study curves $\gamma$ in the punctual Hilbert scheme
$M_n(x) \subset \hil{n}$ such that $\gamma \sim \beta_n$. 
We shall establish some notations, and prove three technical lemmas.
We identify $M_n(x)$ with $\Hilb^n(\Cee^2, 0)$ throughout this section.
Let $R = \Cal O_{\Bbb C^2, 0}$ be the local ring of
$\Cee^2$ at the origin, and $\frak m=(u, v)$ be the maximal ideal of
$R$ where $u, v$ form a coordinate system of $\Cee^2$ at the origin.
Let $\eta \in \Hilb^n(\Bbb C^2, 0)$.
It is known (see \cite{Iar}) that $I_{\eta, 0} \supset \frak m^n$
and there is an embedding 
$\tau\colon\, \Hilb^n(\Bbb C^2, 0) \to \Grass(R/\frak m^n, n)$
where $R/\frak m^n$ is considered as a $\Bbb C$-vector space
of dimension $n+1\choose 2$, and $\tau$ maps
$\eta\in \Hilb^n(\Bbb C^2, 0)$ to the quotient
$R/\frak m^n \to R/I_{\eta, 0}=\Cal O_{\eta, 0} \to 0$.
For simplicity, put $\Bbb G = \Grass(R/\frak m^n, n)$ and
$\ell_1 ={n+1\choose 2}$. Then over $\Bbb G$,
there are a universal quotient bundle $\Cal Q$,
a universal subbundle $\Cal S$ and an exact sequence:
$$\exact{\Cal S}{(R/\frak m^n) \otimes \Cal O_{\Bbb G}}{\Cal Q}. 
\tag 4.1$$
It is known that det$(\Cal Q)$ induces the Pl\"ucker embedding 
$\frak p: \Bbb G \to \Pee((\wedge^{\ell_1 - n} \Cee^{\ell_1})^*).$
So if $\Cal H$ is the hyperplane line bundle over 
$\Pee((\wedge^{\ell_1 -n} \Cee^{\ell_1})^*)$, 
then $\text{det}(\Cal Q) = \frak p^*\Cal H$.

\lemma{4.2} Let $n\ge 2$. 
Fix a point $x \in X$ and a curve $\gamma \subset M_n(x)$.
Regard $\gamma \subset \Hilb^n(\Cee^2, 0)$
via the identification $M_n(x) \cong \Hilb^n(\Cee^2, 0)$. 
Then, $\gamma\sim\beta_n$ if and only if $\gamma \cdot B_n=-2$ 
if and only if $\frak p \circ \tau (\gamma)$ is a line.
\endstatement
\proof
Consider the rank-$n$ vector bundle
$p_{n*}\Cal O_{\Cal Z_n}$ over $\hil{n}$. Following the similar argument
as in the proof of the Claim 2 in the Appendix in \cite{BSG},
we can show that $(p_{n*}\Cal O_{\Cal Z_n})|_{M_n(x)}$ is a quotient of
the trivial bundle $(\Cal O_{M_n(x)})^{\oplus \ell_1}$
where $\ell_1 ={n+1\choose 2}$.
The fiber of $(p_{n*}\Cal O_{\Cal Z_n})|_{M_n(x)}$ over
$\eta \in M_n(x) \subset \hil{n}$ is $H^0(X, \Cal O_\eta) = \Cal O_{\eta,x}$.
Since the fiber of $\tau^*\Cal Q$ over $\eta \in \Hilb^n(\Cee^2, 0)$
is $\Cal O_{\eta, 0}$, we see from  the universality of
the quotient bundle that
$(p_{n*}\Cal O_{\Cal Z_n})|_{M_n(x)} \cong \tau^*\Cal Q$.
It is known from the Appendix in \cite{BSG} that
$c_1(p_{n*}\Cal O_{\Cal Z_n})=-B_n/2$.
Combining with $\text{det}(\Cal Q) = \frak p^*\Cal H$ yields
$$\align
&\gamma \cdot B_n
 = -2 \gamma \cdot c_1(p_{n*}\Cal O_{\Cal Z_n})
 = -2 \gamma \cdot c_1((p_{n*}\Cal O_{\Cal Z_n})|_{M_n(x)})\\
=& -2 \gamma \cdot c_1(\tau^*\Cal Q)
 = -2 \gamma \cdot \tau^*\frak p^*c_1(\Cal H)  
= -2 (\frak p \circ \tau)(\gamma) \cdot c_1(\Cal H). \tag 4.3
\endalign$$

Next, since $\gamma \subset M_n(x)$, 
we have $\gamma \cdot D_C=0=\beta_n\cdot D_C$ for 
every real surface $C$ in $X$. By Lemma 2.10 and the Poincar\'e duality,
$\gamma \sim \beta_n$ if and only if $\gamma \cdot B_n=-2=\beta_n\cdot B_n$.
By (4.3), $\gamma \cdot B_n=-2$ if and only if 
$\frak p \circ \tau (\gamma)$ is a line.
\endproof

\lemma{4.4} Let $k=n(n-1)/2$, and $\Bbb L \subset \Bbb G
= \Grass(R/\frak m^n, n)$ be a curve. Then, $\Bbb L$ is mapped to
a line by $\frak p$ if and only if $\Bbb L$ is of the form
$$\Bbb L = \{ \text{\rm Span}(\overline f_1, \ldots, \overline f_{k-1},
\lambda\overline f_k +\mu \overline f_{k+1})\, |\,
\lambda, \mu \in \Bbb C, |\lambda|+|\mu|\neq0\} \tag 4.5$$
where $f_1, \ldots, f_{k+1}\in R$ are polynomials such that 
$\overline f_1, \ldots, \overline f_{k+1}$ are linearly independent
in the vector space $R/\frak m^n \cong \Bbb C^{n+k}$. 
\endstatement
\proof
 From the definition of $\frak p$, it is clear that if $\Bbb L$
is of the form (4.5), then $\frak p(\Bbb L)$ is a line. 
Conversely, assume that $\frak p(\Bbb L)$ is a line.
We see from $\text{det}(\Cal Q) = \frak p^*\Cal H$ 
that $\det (\Cal Q|_{\Bbb L})\cong \Cal O_{\Bbb L}(1)$.
Restricting (4.1) to $\Bbb L$, we obtain an exact sequence:
$$\exact{\Cal S|_{\Bbb L}}{(R/\frak m^n) \otimes \Cal O_{\Bbb L}}
{\Cal Q|_{\Bbb L}}.  \tag 4.6$$
So $\det (\Cal S|_{\Bbb L} )\cong \Cal O_{\Bbb L}(-1)$. 
 From the injection $\Cal S|_{\Bbb L} \rightarrow 
(R/\frak m^n) \otimes \Cal O_{\Bbb L}$, we obtain 
$\Cal S|_{\Bbb L} \cong \bigoplus_{i=1}^{k}\Cal O_{\Bbb L}(a_i)$ 
with $\sum_{i=1} ^{k}a_i= -1$ and $a_i \le 0$. Thus $\Cal S|_{\Bbb L} 
\cong \Cal O_{\Bbb L}^{\oplus (k-1)}\oplus \Cal O_{\Bbb L}(-1)$. 
Now the injection $\Cal O_{\Bbb L}^{\oplus (k-1)} \rightarrow 
(R/\frak m^n) \otimes \Cal O_{\Bbb L}$ gives rise to 
linearly independent elements 
$\overline f_1, \ldots, \overline f_{k-1}$ in $R/\frak m^n$. 
The injection $\Cal O_{\Bbb L}(-1)\rightarrow  
(R/\frak m^n) \otimes \Cal O_{\Bbb L}$ gives rise to
two linearly independent elements $\overline f_k$ and $\overline f_{k+1}$ 
in $R/\frak m^n$ such that $\overline f_1, \ldots, \overline f_{k-1}, 
\lambda \overline f_k + \mu \overline f_{k+1}$ are
linearly independent when $|\lambda|+|\mu|\neq 0$. 
So ${\Bbb L}$ is of the form (4.5).
\endproof

\lemma{4.7} Let $x \in X$ and $\eta\in M_{n+1}(x)$.
Let $\xi_1, \xi_2 \in M_{n}(x)$ and $\xi_1 \ne \xi_2$.
If $\xi_1, \xi_2 \subset \eta$, then $I_{\eta, x}\supset \frak m_x^n$
where $\frak m_x$ is the maximal ideal of $\Cal O_{X,x}$.
\endstatement
\proof
We have $I_{\eta,x} \subset I_{\xi_1, x} \cap I_{\xi_2,x} \subset
I_{\xi_2,x}$ and $I_{\xi_1, x} \cap I_{\xi_2,x} \ne I_{\xi_2,x}$.
Since the colength of $I_{\eta,x}$ is $n+1$ and 
the colength of $I_{\xi_2,x}$ is $n$, 
we obtain $I_{\eta,x} =I_{\xi_1, x} \cap I_{\xi_2,x}$. 
Since $\frak m_x^n\subset I_{\xi_1, x}$ and $\frak m_x^n\subset I_{\xi_2,x}$,
we get $\frak m_x^n\subset I_{\xi_1,x}\cap I_{\xi_2,x}=I_{\eta,x}$.
\endproof

Now we are able to characterize all the curves $\gamma$ in
the punctual Hilbert scheme $M_n(x)$ which is homologous to $\beta_n$.
By the Proposition 2.1 of \cite{ES2}, the fiber $(\psi_{n+1})^{-1}(\eta, x)$
over $(\eta, x) \in \Cal Z_{n+1} \subset \hil{n+1} \times X$
is isomorphic to $\Pee(\omega_{\eta_x} \otimes \Cee_x)$.

\proposition{4.8} Let $n\ge 2$, and $X$ be simply-connected.
Let $x \in X$, and $\gamma$ be a curve in $M_n(x)$.
Then, $\gamma \sim \beta_n$ if and only if $\gamma=f_{n+1}(C)$ 
where $C$ is a line in the fiber $(\psi_{n+1})^{-1}(\eta, x) \cong 
\Pee(\omega_{\eta} \otimes \Cee_x)$ for some $\eta\in M_{n+1}(x)$.
\endstatement
\proof 
(i) Fix the identification $M_n(x) \cong \Hilb^n(\Cee^2, 0)$ in (2.4). 
Let $C$ be a line in $(\psi_{n+1})^{-1}(\eta, x) \cong 
\Pee(\omega_{\eta} \otimes \Cee_x)$ for some $\eta\in M_{n+1}(x)$.  
By Lemma 4.7, $I_{\eta, x}\supset \frak m^n$.
Choose polynomials $f_1, \ldots, f_{k-1} \in I_{\eta}$ such that 
$\overline f_1, \ldots, \overline f_{k-1}$ form a basis of the 
$(k-1)$-dimensional subspace $I_{\eta, x}/\frak m^n\subset R/\frak m^n$. 
Now, $I_{\eta, x} =
\Cee f_1+\ldots+\Cee f_{k-1}+\frak m^n$. So $I_{\eta} = 
I_{\eta, x} \cap \Cee[u, v]= \Cee f_1+\ldots+\Cee f_{k-1}+(u, v)^n$.

Next, since $\Pee(\omega_{\eta} \otimes \Cee_x) \cong 
\Pee(\Hom(\Cee_x, \Cal O_{\eta})^*)$, the line $C$ in 
$(\psi_{n+1})^{-1}(\eta, x) \cong \Pee(\Hom(\Cee_x, \Cal O_{\eta})^*)$
is generated by two different points $[\epsilon_1], [\epsilon_2]
\in \Pee(\Hom(\Cee_x, \Cal O_{\eta})^*)$ (so that the two morphisms 
$\epsilon_1, \epsilon_2 \in \Hom(\Cee_x, \Cal O_{\eta}) 
= \Hom_R(R/\frak m, \Cal O_{\eta})$ are linearly independent). 
Choose $f_k, f_{k+1} \in R$ such that 
$[f_k]= \epsilon_1(1)$ and $[f_{k+1}] =\epsilon_2(1)$ in $ R/I_{\eta, x}$.
Note that $\overline f_1, \ldots, \overline f_{k+1} 
\in R/\frak m^n$ are linearly independent. 

Let $\xi\in \gamma = f_{n+1}(C)$, and assume that the point 
$(\xi, \eta) \in C$ corresponds to 
$[\lambda \epsilon_1+\mu \epsilon_2] \in 
\Pee(\Hom(\Cee_x, \Cal O_{\eta})^*)$. Then we have the exact sequence:
$$0\longrightarrow \Cee_x\mapright{\lambda \epsilon_1+\mu \epsilon_2}
\Cal O_{\eta}\longrightarrow \Cal O_{\xi}\longrightarrow 0.$$
So $I_{\xi} = I_\eta + \Cee(\lambda f_k+\mu f_{k+1})
= \Cee f_1+\ldots+\Cee f_{k-1}+ \Cee(\lambda f_k+\mu f_{k+1})+(u, v)^n$. 
Thus $\tau(\gamma) = \Bbb L$ which is of the form (4.5). 
By Lemma 4.4 and Lemma 4.2, $\gamma \sim \beta_n$.

(ii) Conversely, assume $\gamma \sim \beta_n$. By Lemma 4.2 and 
Lemma 4.4, $\tau(\gamma) = \Bbb L$ which is of the form (4.5). 
Let $u\colon R/\frak m^n \to R/\frak m^n$ be the map defined by
multiplication with $u$. 
Let $\varphi_u\colon R/\frak m^n\otimes \Cal O_{\Bbb L} \to 
R/\frak m^n\otimes\Cal O_{\Bbb L}$ be the induced morphism,
and consider the restriction of $\varphi_u$ to the subbundle 
$\Cal S|_{\Bbb L} \cong \Cal O_{\Bbb L}(-1)\oplus 
\Cal O_{\Bbb L}^{\oplus (k-1)}$ in (4.6). 
Since $\Bbb L = \tau(\gamma)$, 
$\text{\rm Span}(\overline f_1, \ldots, \overline f_{k-1},
\lambda\overline f_k +\mu \overline f_{k+1})$ is an ideal of $R/\frak m^n$
whenever $|\lambda|+|\mu| \neq 0$. 
So we see from the proof of Lemma 4.4 that $\varphi_u$ maps 
$\Cal O_{\Bbb L}(-1)\oplus \Cal O_{\Bbb L}^{\oplus (k-1)}$ to itself. 
Since $\Hom(\Cal O_{\Bbb L}, \Cal O_{\Bbb L}(-1))=0$, $\varphi_u$ maps
$\Cal O_{\Bbb L}^{\oplus (k-1)}$ to $\Cal O_{\Bbb L}^{\oplus (k-1)}$. 
The same statement holds when we replace $u$ by $v$. 
Hence, $\text{\rm Span}(\overline f_1, \ldots, \overline f_{k-1})$ 
is an ideal of $R/\frak m^n$, and $\Cee f_1+\ldots+\Cee f_{k-1}+\frak m^n$
is an ideal of $R$ with colength $n+1$. Let $\eta \in M_{n+1}(x)$ 
be defined by the ideal $I_\eta {\overset \text{def} \to =} 
\Cee f_1+\ldots+\Cee f_{k-1}+ (u, v)^n$.

Since $\Hom(\Cal O_{\Bbb L}(-1), \Cal O_{\Bbb L}(-1)) 
=\Cee \cdot \Id_{\Cal O_{\Bbb L}(-1)}$, we have 
$\varphi_u|_{\Cal O_{\Bbb L}(-1)} = c \cdot \Id_{\Cal O_{\Bbb L}(-1)}$ 
for some $c \in \Cee$. Then $(\varphi_u|_{\Cal O_{\Bbb L}(-1)})^n
= c^n \cdot \Id_{\Cal O_{\Bbb L}(-1)}$. 
However, $(\varphi_u)^n =0$ since $u^n\in \frak m^n$. 
So $c=0$, and  $\varphi_u$ maps $\Cal O_{\Bbb L}(-1) \oplus 
\Cal O_{\Bbb L}^{\oplus (k-1)}$ to $\Cal O_{\Bbb L}^{\oplus (k-1)}$. 
In particular, $u {\overline f}_k, u {\overline f}_{k+1} \in
\text{\rm Span}(\overline f_1, \ldots, \overline f_{k-1})$. 
Hence $uf_k, uf_{k+1} \in I_{\eta}$. 
Similarly, $vf_k, vf_{k+1} \in I_{\eta}$.

Define two $R$-morphisms $R \to \Cal O_{\eta}$ by sending $1 \in R$ to 
$[f_k], [f_{k+1}]\in \Cal O_{\eta}$ respectively. 
They induce two $R$-morphisms $\epsilon_1, 
\epsilon_2 \colon \Cee_x = R/\frak m \to \Cal O_{\eta}$
since $uf_k$, $uf_{k+1}$, $vf_k$, $vf_{k+1} \in I_{\eta}$. 
Since $\overline f_1, \ldots, \overline f_{k-1}, \overline f_{k},
\overline f_{k+1} \in R/\frak m^n$ are linearly independent, 
so are the two morphisms
$\epsilon_1, \epsilon_2 \in \Hom(\Cee_x, \Cal O_{\eta})$. 
Let $C$ be the line in $(\psi_{n+1})^{-1}(\eta, x)
\cong \Pee(\Hom(\Cee_x, \Cal O_{\eta})^*)$ generated by the points 
$[\epsilon_1], [\epsilon_2] \in \Pee(\Hom(\Cee_x, \Cal O_{\eta})^*)$.
As in the last paragraph in part (i) above, 
the points in $f_{n+1}(C)$ correspond to the ideals $\Cee f_1+\ldots
+\Cee f_{k-1}+\Cee (\lambda f_k +\mu f_{k+1})+(u,v)^n$,
$\lambda, \mu \in \Cee$ with $|\lambda|+|\mu| \neq 0$.
These same ideals define the points in $\gamma$ since $\tau(\gamma) = 
\Bbb L$ is given by (4.5). So $\gamma = f_{n+1}(C)$.
\qed

Next, we study an arbitrary curve $\gamma \subset \hil{n}$ 
homologous to $\beta_n$. 

\lemma{4.9} 
{\rm (i)} A curve $\gamma \subset \hil{n}$ homologous to $\beta_n$
is a smooth rational curve;
\par
{\rm (ii)} A curve $\gamma \subset \hil{n}$ is homologous to $\beta_n$
if and only if $\gamma = \gamma_k + \xi$ where $2 \le k \le n$,
$\gamma_k$ is a curve in the punctual Hilbert scheme $M_k(x)$
for some $x \in X$, $\gamma_k \sim \beta_k$ in $\hil{k}$,
and $\xi \in \hil{n-k}$ is a fixed point with $x \not \in \Supp(\xi)$.
\endstatement
\proof
(i) By Remark 3.3, $(D_C - B_n/2)$ is very ample for some divisor $C$ on
$X$. Since $\gamma \cdot (D_C-B_n/2) =\beta_n \cdot (D_C-B_n/2) = 1$,
$\gamma$ is a smooth rational curve.

(ii) First of all, assume $\gamma \sim \beta_n$.
Since $\gamma$ is irreducible and
$\gamma \cdot D_C = \beta_n \cdot D_C = 0$,
$\rho_n(\gamma)$ is a point in $X^{(n)}$.
Let $\rho_n(\gamma)=n_1x_1+\ldots+n_mx_m$ where $n_i\ge 1$,
$\sum_{i=1}^m n_i = n$, and $x_1, \ldots, x_m$ are distinct points.
Then $\gamma$ is contained in the product
$\prod_{i=1}^m M_{n_i}(x_i)\subset \hil{n}$.
Since $M_{n_i}(x_i) \cong \Hilb^{n_i}(\Bbb C^2, 0)$,
$H_1(M_{n_i}(x_i), \Bbb Z) =0$ by the Theorem 1.1 (i) of \cite{ES1}.
By the K\"unneth formula, $\gamma \sim \sum_{i=1}^m (p_i(\gamma)+\xi_i)$ 
in $\prod _{i=1}^m M_{n_i}(x_i)$,
where $p_i$ is the $i$-th projection of $\prod_{i=1}^m M_{n_i}(x_i)$, 
and $\xi_i \in \hil{n-n_i}$ with
$\Supp(\xi_i)\cap\Supp\big(p_i(\gamma)\big)=\emptyset$. So
$1 =(D_C-B_n/2) \cdot \gamma
=\sum_{i=1}^m(D_C-B_n/2)\cdot (p_i(\gamma)+\xi_i).$
Since $(D_C-B_n/2)$ is very ample, 
there exists exactly one $i$ such that
$p_i(\gamma)$ is a curve in $\hil{n_i}$ and
$p_j(\gamma)$ is a point in $\hil{n_j}$ for all $j \neq i$.
We may assume $i=1$.
Let $k=n_1$ and $x=x_1$. Then $\gamma=\gamma_k+\xi$ where
$\gamma_k \subset M_k(x)$ and $\xi\in \hil{n-k}$ is fixed
with $x \notin \Supp(\xi)$.
It remains to prove $\gamma_k \sim \beta_k$ in $\hil{k}$.

\medskip\noindent
{\bf Claim.} {\it Let $\w x \in X$ be a fixed point, 
and $\alpha \subset \hil{n-1}$ be 
a curve such that $\w x \notin \Supp(\xi)$ for every $\xi\in \alpha$. 
Then, $(\alpha+\w x)\cdot B_n=\alpha\cdot B_{n-1}$.}
\proof 
Define
$\W \alpha=\{\,(\xi, \xi+\w x)\,|\, \xi\in \alpha\}
\subset \hil{n-1, n}.$
Since $\w x \notin \Supp(\alpha)$,
the restriction $g_n|_{\W \alpha}: \W \alpha \to \alpha+\w x$
is an isomorphism. So $(g_n)_*\W\alpha=\alpha+\w x$.
Similarly, $(f_n)_*\W \alpha=\alpha$.
Let $E_n =\{ (\xi, \eta)\in \hil{n-1, n}\, |\,\Supp(\xi)
=\Supp(\eta), \xi\subset \eta\}$.
Since $\w x \notin \Supp(\alpha)$, we obtain $\W \alpha\cap E_n=\emptyset$.
Since $((g_n)^*B_n-(f_n)^*B_{n-1})=2E_n$ (see \cite{Leh}), 
we have $\W\alpha\cdot ((g_n)^*B_n-(f_n)^*B_{n-1})=0$. 
Therefore, we obtain
$(\alpha+\w x)\cdot B_n=(g_n)_*\W\alpha\cdot B_n
=\W\alpha\cdot (g_n)^* B_n=\W \alpha\cdot (f_n)^*B_{n-1}
=(f_n)_* \W\alpha\cdot B_{n-1}=\alpha\cdot B_{n-1}$.
\endproof

By deforming $\xi$ in $\hil{n-k}$, we can assume that 
$\xi$ is $\w x_1+\ldots+\w x_{n-k}$
where $\w x_1,\ldots,\w x_{n-k}$ are distinct points different from $x$.
From the Claim, we get  
$$(\gamma_k+\w x_1+\ldots+\w x_{n-k})\cdot B_n
=(\gamma_k+\w x_1+\dots+\w x_{n-k-1})\cdot B_{n-1}
=\ldots=\gamma_k\cdot B_k.  \tag 4.10$$
So $\gamma_k\cdot B_k = \gamma\cdot B_n = \beta_n \cdot B_n= -2$.
Since $\gamma_k \subset M_k(x)$, $\gamma_k \sim \beta_k$ by Lemma 4.2.

Conversely, suppose $\gamma = \gamma_k + \xi$ where $2 \le k \le n$,
$\gamma_k \subset M_k(x)$ for some $x \in X$, 
$\gamma_k \sim \beta_k$ in $\hil{k}$,
and $\xi \in \hil{n-k}$ is fixed with $x \not \in \Supp(\xi)$.
We may assume that $\xi$ consists of $(n-k)$ distinct points.
By (4.10) and (3.13), $\gamma \cdot B_n = \gamma_k \cdot B_k
= \beta_k \cdot B_k =-2$. Then the same argument as
in the proof of Lemma 4.2 yields $\gamma \sim \beta_n$.
\endproof

\theorem{4.11} Let $n\ge 2$, and $X$ be a simply-connected smooth 
projective surface. Then, a curve $\gamma$ in $\hil{n}$ is homologous to 
$\beta_n$ if and only if $\gamma=f_{n+1}(C)$ where $C$ is 
a line in the projective space $(\psi_{n+1})^{-1}(\eta, x)$ 
for some $(\eta, x) \in \Cal Z_{n+1}$. Moreover, in this case, 
$(\eta, x)$ and $C$ are uniquely determined by $\gamma$.
\endstatement
\proof
(i) First, assume $\gamma \sim \beta_n$. By Lemma 4.9 (ii),
we know that $\gamma=\gamma_k+\xi$
where $\gamma_k \subset M_k(x)$ for some $x \in X$,
$\xi\in \hil{n-k}$ is fixed with $x \not\in \Supp(\xi)$,
and $\gamma_k \sim \beta_k$ in $\hil{k}$. By Proposition 4.8,
$\gamma_k = f_{k+1}(\W C)$ where $\W C$ is a line in 
$(\psi_{k+1})^{-1}(\eta_{k+1}, x)$ for some $\eta_{k+1} \in M_{k+1}(x)$.
Since $\Supp(\eta_{k+1}) = \{x\} \not\subset \Supp(\xi)$, we choose
$\eta {\overset \text{def} \to =} \eta_{k+1}+\xi \in \hil{n+1}$
and $(\eta, x)\in \Cal Z_{n+1}$.
By the definition of $\psi_{n+1}$ in Notation 2.5,
$$\align
     &\,\, (\psi_{n+1})^{-1}(\eta, x)
      =\{ (\zeta, \eta)\in \hil{n, n+1} \, |\, \zeta\subset \eta
      \text{ and } \Supp(I_\zeta/I_\eta) = x \}   \\
=    &\,\, \{ (\zeta_{k} + \xi, \eta_{k+1} + \xi) \in \hil{n, n+1} \, |\,
      \zeta_{k} \subset \eta_{k+1} \}   \\
\cong&\,\, \{ (\zeta_{k}, \eta_{k+1}) \in \hil{k, k+1} \, |\,
      \zeta_{k} \subset \eta_{k+1} \}    
      =(\psi_{k+1})^{-1}(\eta_{k+1}, x). 
\endalign$$
Put $C = \{ (\zeta_k + \xi, \eta_{k+1}+\xi)| (\zeta_k, \eta_{k+1}) \in \W C \}
= \{ (\zeta_k + \xi, \eta)| (\zeta_k, \eta_{k+1}) \in \W C \}$, i.e., $C$
is the line in $(\psi_{n+1})^{-1}(\eta, x)$ corresponding to the line $\W
C$ in
$(\psi_{k+1})^{-1}(\eta_{k+1}, x)$ via the above isomorphism
$(\psi_{n+1})^{-1}(\eta, x) \cong (\psi_{k+1})^{-1}(\eta_{k+1}, x)$.
Then, we see that $\gamma=\gamma_k+\xi=f_{k+1}(\W C)+\xi=f_{n+1}(C)$.

Conversely, suppose $\gamma=f_{n+1}(C)$ where $C$ is a line
in $(\psi_{n+1})^{-1}(\eta, x)$ 
for some $(\eta, x) \in \Cal Z_{n+1}$. Assume $\ell(\eta_x)=k+1$.
Let $\eta_{k+1} = \eta_x \in \hil{k+1}$.
Reversing the argument in the preceding paragraph and
using Lemma 4.9 (ii), we get $\gamma \sim \beta_n$.

(ii) To show the uniqueness of $(\eta, x)$ and $C$, 
let $\xi_1, \xi_2 \in \gamma$ with $\xi_1 \ne \xi_2$. 
Then, $\xi_1, \xi_2 \subset \eta$.
Since $\ell(\eta) = 1+\ell(\xi_1) = 1+\ell(\xi_2)$,
we get $I_\eta = I_{\xi_1} \cap I_{\xi_2}$. 
So $\eta$ and hence the point $x$ are uniquely determined by $\gamma$.
Since $f_{n+1}|_{(\psi_{n+1})^{-1}(\eta, x)}: (\psi_{n+1})^{-1}(\eta, x) 
\to f_{n+1}((\psi_{n+1})^{-1}(\eta, x))$ is an isomorphism,
$C$ is also uniquely determined by $\gamma$.
\endproof

Our next goal is to give a global description of the moduli space
$\frak M(\beta_n)$ of all the curves in $\hil{n}$ homologous to $\beta_n$.
Note that $\frak M(\beta_n)$ is the union of certain irreducible
components in the Hilbert scheme of curves in $\hil{n}$. We shall use
the concept of the relative Grassmannian scheme (see \cite{Sim}). 

\theorem{4.12} Let $n\ge 2$, $X$ be a simply-connected 
smooth projective surface, and 
$\omega_{\Cal Z_{n+1}}$ be the dualizing sheaf of 
$\Cal Z_{n+1} \subset X^{[n+1]} \times X$.
Then, there exists a bijective morphism from the relative Grassmannian
$\text{\rm Grass}(\omega_{\Cal Z_{n+1}}, 2)$ over $\Cal Z_{n+1}$
to the moduli space $\frak M(\beta_n)$ of
all the curves in $\hil{n}$ homologous to $\beta_n$.
\endstatement
\proof
Put $\W {\Bbb G} = \text{Grass}(\omega_{\Cal Z_{n+1}},2)$ for simplicity.
Define $\Gamma = \{ ([C], (\xi, \eta)) \in 
\W {\Bbb G} \times \hil{n, n+1}|\,\,
C \text{ is a line in } (\psi_{n+1})^{-1}(\eta, x)
\text{ and } (\xi, \eta) \in C \text{ with } 
\Supp(I_\xi/I_\eta) = \{x\} \}$. By the Proposition 2.1 of \cite{ES2}, 
$\Pee(\omega_{\Cal Z_{n+1}}) \cong \hil{n, n+1}$ 
over $\Cal Z_{n+1}$. 
So $\Gamma$ is the universal subscheme in
$\W {\Bbb G} \times_{\Cal Z_{n+1}}
\Pee(\omega_{\Cal Z_{n+1}})$ which is flat over $\W {\Bbb G}$. 
Consider 
$$\alpha: \W {\Bbb G} \times_{\Cal Z_{n+1}} \hil{n, n+1} \subset
\W {\Bbb G} \times \hil{n, n+1} \text{ }
\overset {\text{Id}_{\W {\Bbb G}} \times f_{n+1}} \to \longrightarrow
\text{ } \W {\Bbb G} \times \hil{n}.$$
Note that the fiber of $\W {\Bbb G}$ over $(\eta, x)\in \Cal Z_{n+1}$ is
$\W {\Bbb G}|_{(\eta, x)}=\text{Grass}(\omega_{\eta_x} \otimes \Cee_x,2)$. So
$$\align
 \Gamma|_{(\eta, x)}
&\subset \W {\Bbb G}|_{(\eta, x)} \times
  \Pee(\omega_{\Cal Z_{n+1}})|_{(\eta, x)}   
= \text{Grass}(\omega_{\eta_x} \otimes \Cee_x,2)
  \times \Pee(\omega_{\eta_x} \otimes \Cee_x)   \\
&= \text{Grass}(\omega_{\eta_x} \otimes \Cee_x,2) \times
  (\psi_{n+1})^{-1}(\eta, x).  \tag 4.13
\endalign$$
Since $f_{n+1}$ maps $(\psi_{n+1})^{-1}(\eta, x)$
isomorphically to $f_{n+1}((\psi_{n+1})^{-1}(\eta, x))$, we see that
$\alpha|_{\Gamma}: \Gamma \to \alpha(\Gamma)$ is bijective.
Moreover, assuming that $q \in \W {\Bbb G}$ is mapped to $(\eta, x)$
via the natural projection $\W {\Bbb G} \to \Cal Z_{n+1}$
and $q$ corresponds to the line
$$\{q\} \times C \overset \text{def} \to = \Gamma|_q
\subset \{q\} \times \Pee(\omega_{\eta_x} \otimes \Cee_x)
= \{q\} \times (\psi_{n+1})^{-1}(\eta, x)$$
via the universal subscheme $\Gamma \subset
\W {\Bbb G} \times_{\Cal Z_{n+1}} \Pee(\omega_{\Cal Z_{n+1}})$,
then we see from (4.13) that the same point $q \in \W {\Bbb G}$
also corresponds to the curve
$$\alpha(\Gamma)|_q = \alpha(\Gamma|_q) = \alpha(\{q\}\! \times\! C)
= \{q\}\! \times \!f_{n+1}(C) \subset
\{q\} \times f_{n+1}((\psi_{n+1})^{-1}(\eta, x))$$
via the subscheme $\alpha(\Gamma) \subset \W {\Bbb G} \times \hil{n}$.
By Theorem 4.11, we conclude that via the subscheme
$\alpha(\Gamma) \subset \W {\Bbb G} \times \hil{n}$ which is flat
over $\W {\Bbb G}$, the relative Grassmannian $\W {\Bbb G}$ parameterizes
curves in $\hil{n}$ homologous to $\beta_n$. By the universal property of
the Hilbert scheme of curves, we obtain a morphism
$\Lambda: \W {\Bbb G} \to \frak M(\beta_n)$. Moreover,
we see from Theorem 4.11 that the morphism $\Lambda$ is bijective.
\endproof

Finally, we study a top stratum in $\text{Grass}(\omega_{\Cal Z_{n+1}},2)$.
Let $x\in X$ and $\xi = x_1+\ldots+x_{n-2} \in \hil{n-2}$ be fixed such that
$x, x_1, \ldots, x_{n-2}$ are distinct points.
Then the curve $M_2(x)+\xi$ is homologous to $\beta_n$. 
The set of all such curves has dimension
$$\# (\hbox{moduli of }x)+
\# (\hbox{moduli of }\xi)=2+2(n-2)=2n-2.$$
We shall prove that this set is the top stratum in 
$\text{Grass}(\omega_{\Cal Z_{n+1}}, 2)$. For this purpose, 
we now give a more detailed geometric description of this set as follows.
Note that $\frak m_x^2$ is the only colength-$3$ ideal
in $\Cal O_{X, x}$ whose minimal number of generators is $3$. 
Let $\theta(x)$ be the unique length-$3$ $0$-dimensional subscheme of $X$ 
corresponding to $\frak m_x^2$.
Let $\Delta_{n-1}=\{ (x_0, x_1, \ldots, x_{n-2}) \in X^{n-1}|
x_0, x_1, \ldots, x_{n-2} \text{ are not distinct} \} $.
Define $\Psi: X^{n-1}-\Delta_{n-1} \to \hil{n+1}$
by $\Psi(x, x_1, \ldots, x_{n-2}) = \theta(x)+\sum_{i=1}^{n-2} x_i$.
Let $\theta(x)+\xi \in \im(\Psi)$ where $\xi = \sum_{i=1}^{n-2} x_i$.
Then, $(\pi_{n+1})^{-1}(\theta(x)+\xi)$ consists of
$(n-1)$ distinct points in $\Cal Z_{n+1}$:
$(\theta(x)+\xi, x)$, $(\theta(x)+\xi, x_i)$, $1 \le i \le n-2$.
Note that 
$\text{Grass}(\omega_{\Cal Z_{n+1}}, 2)|_{(\theta(x)+\xi, x)}
=\text{Grass}(\omega_{\theta(x)} \otimes \Cee_x, 2)
\cong \Grass(\Bbb C^2, 2) = pt$,
where we have used (3.5) in \cite{ES2}. In fact,
$\text{Grass}(\omega_{\Cal Z_{n+1}}, 2)|_{(\theta(x)+\xi, x)}$
is the point in 
$\text{Grass}(\omega_{\Cal Z_{n+1}}, 2)$
corresponding to the curve $M_2(x) + \xi$ in $\hil{n}$.
Also, $\text{Grass}(\omega_{\Cal Z_{n+1}}, 2)|_{(\theta(x)+\xi, x_i)}
\cong \Grass(\Bbb C, 2)= \emptyset$. So 
$\text{Grass}(\omega_{\Cal Z_{n+1}}, 2)
|_{(\pi_{n+1})^{-1}(\im(\Psi))}  $
is exactly the subset of $\text{Grass}(\omega_{\Cal Z_{n+1}}, 2)$
described at the beginning of this paragraph. 

\proposition{4.14} {\rm (i)} The subset 
$\text{\rm Grass}(\omega_{\Cal Z_{n+1}}, 2)
|_{(\pi_{n+1})^{-1}(\im(\Psi))}$ has dimension equal to $(2n-2)$, and
is a top stratum inside $\text{\rm Grass}(\omega_{\Cal Z_{n+1}}, 2)$;

{\rm (ii)} Let $\gamma = M_2(x)+ x_1+\ldots+x_{n-2}$
where $x, x_1, \ldots, x_{n-2}$ are distinct points in $X$.
Then the normal bundle $N_{\gamma\subset \hil{n}}$ of $\gamma$
in $\hil{n}$ is isomorphic to $\Cal O_{\gamma}^{\oplus (2n-2)}
\oplus \Cal O_{\gamma}(-2)$.
\endstatement
\proof  
(i) We need to show that the complement of 
$\text{\rm Grass}(\omega_{\Cal Z_{n+1}}, 2)|_{(\pi_{n+1})^{-1}(\im(\Psi))}$
has dimension less than $(2n-2)$.
Let $(\xi, x) \in \Cal Z_{n+1} - (\pi_{n+1})^{-1}(\im(\Psi))$ such that 
$\Grass(\omega_{\Cal Z_{n+1}}, 2)|_{(\xi, x)} \ne \emptyset$.
Then, $\ell(\xi_x)\ge 3$; moreover, when $\ell(\xi_x) = 3$,
we must have $\xi_x = \theta(x)$.
Since $\xi = \pi_{n+1}(\xi, x) \not\in \im(\Psi)$, 
either $\ell \,\,{\overset \text{def} \to =} \,\, \ell(\xi_x)\ge 4$ 
or $\xi_x = \theta(x)$ and $\xi= \theta(x)+\xi'$ 
where $x\notin \Supp(\xi')$, $\xi' \in \hil{n-2}$ and 
$|\Supp(\xi')| \le (n-3)$. Let $S_1$ and $S_2$ consist of all 
$(\xi, x) \in \Cal Z_{n+1} - (\pi_{n+1})^{-1}(\im(\Psi))$
satisfying these two types of conditions respectively. 
Now the dimension of $\text{\rm Grass}(\omega_{\Cal Z_{n+1}}, 2)|_{S_2}$
is at most 
$$\# (\hbox{moduli of }x)+
\# (\hbox{moduli of }\xi')=2+[2(n-3)+1] < 2n-2$$
since $\text{Grass}(\omega_{\Cal Z_{n+1}}, 2)|_{(\theta(x)+\xi', x)}$
is a point when $x \not \in \Supp(\xi')$. To estimate the dimension of 
$\text{\rm Grass}(\omega_{\Cal Z_{n+1}}, 2)|_{S_1}$, 
let $i=\dim_{\Bbb C}\Hom_X(\Bbb C_x, \Cal O_{\xi_x})$. 
Then, $i=\dim_{\Bbb C} (\omega_{\xi} \otimes \Cee_x)$ by the formula
(3.5) in \cite{ES2} and the Lemma 2 in \cite{E-L}. 
Using the Lemma 7 in \cite{E-L}, we conclude that the dimension of 
$\text{\rm Grass}(\omega_{\Cal Z_{n+1}}, 2)|_{S_1}$ is at most
$$\align
   &\dim (S_1) + \dim \text{Grass}(\omega_{\xi} \otimes \Cee_x, 2)  \\
\le&\left [2+\big(\ell-1-2(i-1) \big) +2(n+1-\ell) \right ]+2(i-2)  \\
=&2n-\ell+1<2n-2
\endalign$$
since $\ell = \ell(\xi_x)\ge 4$ for $(\xi, x) \in S_1$. 
Therefore, the complement of the subset 
$\text{\rm Grass}(\omega_{\Cal Z_{n+1}}, 2)|_{(\pi_{n+1})^{-1}(\im(\Psi))}$
has dimension less than $(2n-2)$.

(ii) In an analytic open neighborhood of $\gamma$, 
$\hil{n}$ splits into $\hil{2} \times \hil{n-2}$.
So $N_{\gamma\subset \hil{n}} \cong N \oplus \Cal O_{\gamma}^{\oplus (2n-4)}$
where $N$ is isomorphic to the normal bundle $\W N$ of $M_2(x)$ in $\hil{2}$
via $\gamma \cong M_2(x)$. Since $\W N \cong \Cal O_{M_2(x)}^{\oplus 2}
\oplus \Cal O_{M_2(x)}(-2)$, our result follows. 
\endproof

By Proposition 4.14 (i) and Theorem 4.12, the dimension of 
$\frak M(\beta_n)$ is $(2n-2)$, one more than the expected dimension 
$c_1(T_{\hil{n}})\cdot \beta_n +\dim\, \hil{n} - \dim\, \text{Aut}(\Bbb P^1) 
=(2n-3)$. Therefore the moduli space of curves in $\hil{n}$ 
homologous to $\beta_n$ is obstructed.

\section{5. Curves homologous to 
$\beta_{\ell}-(n-1)\beta_n$ when $X = \Pee^2$}

In this section, we take $X = \Pee^2$.
We shall characterize all the curves in $\hil{n}$ homologous to
the class $\beta_{\ell}-(n-1)\beta_n$. In addition, we shall compute 
the normal bundles of these curves in $\hil{n}$, 
and describe their moduli space.

First of all, we describe some curves in $\hil{n}$ homologous to
$\beta_{\ell}-(n-1)\beta_n$.

\lemma{5.1} Let $C$ be a line in $X=\Bbb P^2$, and 
consider $C^{(n)} = \Hilb^n(C) \subset \hil{n}$. 

{\rm (i)} Every line in $C^{(n)} \cong \Pee^n$ 
is homologous to $\beta_{\ell}-(n-1)\beta_n$;

{\rm (ii)} $\Cal O_{\hil{n}}(D_{\ell})|_{C^{(n)}} =\Cal O_{C^{(n)}}(1)$.
\endstatement
\proof 
(i) Let $x_1, \ldots, x_{n-1} \in C$ be distinct,
and put $\gamma = C + (x_1+ \ldots + x_{n-1})$. Then 
$\gamma$ is a line in the projective space $C^{(n)} \cong \Pee^n$.
By Lemma 3.11, $\gamma \sim \beta_{\ell}-(n-1)\beta_n$.
So every line in $C^{(n)} \cong \Pee^n$ 
is homologous to $\beta_{\ell}-(n-1)\beta_n$. 
 
(ii) Follows from $\gamma \cdot D_{\ell}|_{C^{(n)}}
= \gamma \cdot D_{\ell} 
= (\beta_{\ell}-(n-1)\beta_n) \cdot D_{\ell} = 1$.
\endproof

Next we shall prove the converse of Lemma 5.1.
Let $\gamma$ be a curve in $\hil{n}$ homologous to
$\beta_{\ell}-(n-1)\beta_n$. Then, $\gamma \cdot (n D_\ell- B_n/2)
= (\beta_{\ell}-(n-1)\beta_n) \cdot (n D_\ell- B_n/2) = 1.$
By Lemma 3.8, $(n D_\ell- B_n/2)$ is very ample.
So $\gamma$ is a smooth rational curve.
Let $p_1$ be the projection from $X^n$ to the first factor. 
Let $C_{\gamma}$ be the union of all the curves in
$p_1\big(\upsilon_n^{-1}\big(\rho_n(\gamma)\big)\big)$ where $\rho_n$ and
$\upsilon_n$ are defined in (2.2). 

\lemma{5.2} With the notations as above,
$C_{\gamma}$ is a line in $X=\Bbb P^2$.
\endstatement
\proof
Suppose $\deg\, (C_{\gamma})\ge 2$. 
Take a point $\xi\in \gamma$ and a smooth point $x\in C_{\gamma}$ such that
$x \notin \Supp(\xi)$. Since $x\in C_{\gamma}\subset
p_1\big(\upsilon_n^{-1}\big(\rho_n(\gamma)\big)\big)$,
there exists $\xi_x \in \gamma$ such that $\rho_n(\xi_x)=n_xx+\eta_x$
where $n_x \ge 1$, $\eta_x\in X^{(n-n_x)}$, and $x\notin \Supp(\eta_x)$.
Choose a line $\ell$ in $X$ missing $\Supp(\eta_x)\cup
\Supp(\xi)$, passing $x$, 
and intersecting $C_{\gamma}$ transversally.
Then the intersection $\ell\cap C_{\gamma}$ is a finite set. 
Since $\deg (C_{\gamma})\ge 2$, $\ell\cap C_{\gamma}$ contains
one more point $y \ne x$.
Hence there exists $\xi_y\in \gamma$ with $y\in \Supp(\xi_y)$. 
Thus $\xi_x, \xi_y \in \gamma\cap D_{\ell}$. Since $y \ne x$, $y \in
\ell$ and $\ell$ misses $\Supp(\eta_x)$, we have
$y\notin \{ x \} \cup \Supp(\eta_x)= \Supp(\xi_x)$. So $\xi_x\neq \xi_y$.
Since $\Supp(\xi) \cap \ell=\emptyset$, we have $\xi \not \in D_{\ell}$.
Since $\xi \in \gamma$, ${\gamma}$ is not contained in $D_{\ell}$.
So $\gamma \cap D_{\ell}$ is a finite set of points.
Since $\xi_x, \xi_y \in \gamma\cap D_{\ell}$ and $\xi_x\neq \xi_y$,
we obtain $\gamma\cdot D_{\ell}\ge 2$, a contradiction to $\gamma\cdot
D_{\ell}=1$.
\endproof

\lemma{5.3} With the notations as above, let $p_{\gamma}$ be a linear
homogeneous polynomial defining the line $C_{\gamma} \subset X=\Pee^2$.
Then, any degree-$(n-1)$ homogeneous polynomial divisible by
$p_{\gamma}$ must vanish at $\xi$ for every $\xi \in \gamma$.
\endstatement
\proof
Let $V_{p_{\gamma}}$ be the subspace of $H^0\big(X, \Cal O_X(n-1)\big)$
consisting of all the degree-$(n-1)$ homogeneous polynomials divisible
by $p_{\gamma}$. Recall from (3.5) that $\varphi_1$
sends $\xi \in \hil{n}$ to the $n(n-1)/2$-dimensional subspace
$H^0(X, \Cal O_X(n-1) \otimes I_{\xi})$ of $H^0\big(X, \Cal O_X(n-1)\big)$.
Since $\gamma\cdot\big((n-1)D_{\ell}- B_n/2\big)
=(\beta_{\ell}-(n-1)\beta_n)\cdot\big((n-1)D_{\ell}- B_n/2\big)=0$,
$\gamma$ is contracted by $\varphi_1$.
So the subspaces $H^0(X, \Cal O_X(n-1)\otimes I_{\xi})\subset
H^0\big(X, \Cal O_X(n-1)\big)$ are independent of $\xi\in \gamma$. Let
$V_{\gamma} =H^0(X, \Cal O_X(n-1)\otimes I_{\xi})$ where $\xi\in \gamma.$
It suffices to prove that $V_{\gamma} =V_{p_{\gamma}}$
as subspaces of $H^0\big(X, \Cal O_X(n-1)\big)$. 

Let $s\in V_{\gamma}$. Then $s$ vanishes at every $\xi\in\gamma$ 
by definition. Since $C_{\gamma}\subset
p_1\big(\upsilon_n^{-1}\big(\rho_n(\gamma) \big)\big)$, 
$s$ vanishes at every point in $C_{\gamma}$.
So the polynomial $s$ is divisible by $p_{\gamma}$. Thus, 
$s \in V_{p_{\gamma}}$, and $V_{\gamma} \subset V_{p_{\gamma}}$.
Since $\dim V_{p_{\gamma}} = h^0(X, \Cal O_X(n-2)) = n(n-1)/2
=h^0(X, \Cal O_X(n-1)\otimes I_{\xi})=\dim V_{\gamma}$
where $\xi \in \gamma$, we have $V_{\gamma} = V_{p_{\gamma}}$.
\endproof

\theorem{5.4} Let $n \ge 2$, and $\ell$ be a line in $X$.
Then, a curve $\gamma\subset\hil{n}$ is homologous to
$\beta_{\ell}-(n-1)\beta_n$ if and only if there is a line $C$ in 
$X$ such that $\gamma$ is a line in $\Hilb^n(C)\subset \hil{n}$.
Moreover, $C$ is uniquely determined by the curve $\gamma$.
\endstatement
\proof
By Lemma 5.1 (i), it remains to prove  the ``only if" part.

Let $C=C_{\gamma}$ where $C_{\gamma}$ is from Lemma 5.2.
First, we claim $\gamma \subset \Hilb^n(C)$. Indeed, let $\xi \in \gamma$.
With the same notations as those in the proof of Lemma 5.3,
we choose $s=p_{\gamma} \cdot q \in V_{p_{\gamma}}$ where
$q\in H^0\big(X, \Cal O_X(n-2)\big)$ and
$q$ doesn't vanish at any point in $\Supp(\xi)$.
By Lemma 5.3, $s$ vanishes at $\xi$.
So
$p_{\gamma}$ must vanish at $\xi$. Therefore, $\xi \subset C_{\gamma} =C$. 
Hence, $\xi\in \Hilb^n(C)=C^{(n)}$ for every $\xi\in \gamma$.
So $\gamma \subset \Hilb^n(C)\subset \hil{n}$.
Next, $\Cal O_{\hil{n}}(D_{\ell})|_{\Hilb^n(C)} =\Cal O_{\Hilb^n(C)}(1)$
by Lemma 5.1~(ii). So viewing $\gamma$ as a curve in $\Hilb^n(C)$, 
we obtain $\gamma \cdot c_1(\Cal O_{\Hilb^n(C)}(1))
= \gamma \cdot D_{\ell}=1$.
Thus $\gamma$ is a line in $\Hilb^n(C)$.

Finally, the uniqueness of the line $C$ follows from the observation
that if $\xi \in \hil{n}$ and $n \ge 2$, then $\xi$ is contained in
at most one line in $X=\Bbb P^2$.
\endproof

In the following, we shall study the normal bundle of a curve $\gamma$
in $\hil{n}$ homologous to $\beta_{\ell}-(n-1)\beta_n$.
By Theorem 5.4, there exists a line $C$ in $X=\Pee^2$ such that
$\gamma$ is a line in $C^{(n)}=\Hilb^n(C)\subset\hil{n}$. In particular,
$N_{\gamma\subset C^{(n)}} \cong \Cal O_{\gamma}(1)^{\oplus (n-1)}$.
So we have the following exact sequence of normal bundles:
$$\exact{\Cal O_{\gamma}(1)^{\oplus (n-1)}}{N_{\gamma\subset \hil{n}}}
{N_{C^{(n)}\subset\hil{n}}|_{\gamma}}.  \tag 5.5$$

\proposition{5.6} Let $n \ge 2$, and $C$ and $\ell$ be lines in 
$X = \Pee^2$. Let $\gamma \subset \hil{n}$ be a curve homologous to
$\beta_\ell - (n-1) \beta_n$. Then,

{\rm (i)} $N_{C^{(n)}\subset\hil{n}}
\cong \pi_{n*}\big(q_n^*\Cal O_X(1)|_{\Cal Z_n} \big)|_{C^{(n)}} \cong
\Cal O_{C^{(n)}}^{\oplus 2}\oplus \Cal O_{C^{(n)}}(-1)^{\oplus (n-2)}$;

{\rm (ii)} $N_{\gamma\subset \hil{n}}=\Cal O_{\gamma}(1)^{\oplus (n-1)}\oplus 
\Cal O_{\gamma}^{\oplus 2}\oplus \Cal O_{\gamma}(-1)^{\oplus (n-2)}$;

{\rm (iii)} the moduli space $\frak M(\beta_\ell - (n-1) \beta_n)$
of all the curves in $\hil{n}$ homologous to $(\beta_\ell - (n-1) \beta_n)$
is unobstructed, i.e., is smooth with the expected dimension.
\endstatement
\proof
(i) By the results in \cite{AIK}, $\Hilb^n(C)$ is the locus of zeros of
a section $s$ of the rank-$n$ locally free sheaf
$\pi_{n*}\big(q_n^*\Cal O_X(1)|_{\Cal Z_n}\big)$ over $\hil{n}$.
Furthermore the section $s$ is transversally regular.
It follows that $N_{C^{(n)}\subset\hil{n}}
\cong \pi_{n*}\big(q_n^*\Cal O_X(1)|_{\Cal Z_n}\big)|_{C^{(n)}}$.

To prove the second isomorphism, 
let $\W {\Cal Z}_n$ be the universal subscheme
in $C^{(n)} \times C \cong \Pee^n \times \Pee^1$. Then, we have 
$\pi_{n*}\big(q_n^*\Cal O_X(1)|_{\Cal Z_n} \big)|_{C^{(n)}} \cong
\w \pi_{n*}\big(\w q_n^*\Cal O_{\Pee^1}(1)|_{\W {\Cal Z}_n}\big)$
where $\w \pi_{n*}: \W {\Cal Z}_n \to \Pee^n$ and 
$\w q_n: \Pee^n \times \Pee^1 \to \Pee^1$ are the projections.
It is known that $\W {\Cal Z}_n \subset \Pee^n \times \Pee^1$ 
is defined by the equation
$a_0U^n + a_1 U^{n-1}V + \ldots + a_n V^n = 0$
where $a_0, a_1, \ldots, a_n$ and $U, V$ are the homogeneous coordinates
on $\Pee^n$ and $\Pee^1$ respectively. So the line bundle
$\Cal O_{\Pee^n \times \Pee^1}(\W {\Cal Z}_n)$ is of type $(1, n)$
in $\Pic(\Pee^n \times \Pee^1) \cong \Zee \oplus \Zee$. 
Let $\w p_n: \Pee^n \times \Pee^1 \to \Pee^n$ be the projection.
Applying $\w p_{n*}$ to the exact sequence
$$0 \to \w q_n^*\Cal O_{\Pee^1}(1) \otimes 
\Cal O_{\Pee^n \times \Pee^1}(-\W {\Cal Z}_n) \to
\w q_n^*\Cal O_{\Pee^1}(1) \to 
\w q_n^*\Cal O_{\Pee^1}(1)|_{\W {\Cal Z}_n} \to 0,$$
we obtain $0 \to \Cal O_{C^{(n)}}^{\oplus 2} \to 
\w \pi_{n*}\big(\w q_n^*\Cal O_{\Pee^1}(1)|_{\W {\Cal Z}_n}\big)
\to \Cal O_{C^{(n)}}(-1)^{\oplus (n-2)} \to 0.$
Since this exact sequence splits, 
we conclude the second isomorphism.

(ii) Follows immediately from (5.5) and the isomorphisms in (i).

(iii) By (ii), $H^1(\gamma, N_{\gamma\subset \hil{n}})=0$
for any curve $\gamma$ homologous to
$(\beta_{\ell}-(n-1)\beta_n)$. By deformation theory, 
$\frak M(\beta_\ell - (n-1) \beta_n)$ is unobstructed.
\endproof  
        
\theorem{5.7} Let $n \ge 2, X=\Pee^2$, and $\ell$ be a line in $X$.
Then the moduli space $\frak M(\beta_{\ell}-(n-1)\beta_n)$ of curves in
$\hil{n}$ homologous to $(\beta_{\ell}-(n-1)\beta_n)$
is isomorphic to the Grassmannian bundle
$\Grass(\Sym^n((T_{(\Bbb P^2)^*})^*), 2)$ over $(\Bbb P^2)^*$.
\endstatement
\proof
First of all, we give a global description to the union $W_n$ of
all the subsets $\Hilb^n(C)\subset \hil{n}$ with
$C$ being lines in $X=\Pee^2$. Consider the incidence variety:
$F=\{ (x, [C])\in \Bbb P^2\times (\Bbb P^2)^*\, |\, x\in C \}
\subset \Bbb P^2\times (\Bbb P^2)^*$. 
We know that $F \cong \Bbb P((T_{(\Bbb P^2)^*})^*)$.
Then we have  $
\Hilb^n(\Bbb P((T_{(\Bbb P^2)^*})^*)/(\Bbb P^2)^*)\cong\Hilb^n(F/(\Bbb P^2)^*) \subset
\Hilb^n(\Bbb P^2\times (\Bbb P^2)^*/(\Bbb P^2)^*)\cong
\hil{n} \times (\Bbb P^2)^*$
of relative Hilbert schemes over $(\Bbb P^2)^*$.
Let $\pi\colon \Hilb^n(\Bbb P((T_{(\Bbb P^2)^*})^*)/(\Bbb P^2)^*)\to (\Bbb P^2)^*$ be 
the projection.
Consider the map
$\alpha\colon \Hilb^n(\Bbb P((T_{(\Bbb P^2)^*})^*)/(\Bbb P^2)^*)
\hookrightarrow \hil{n} \times (\Bbb P^2)^* \to \hil{n}$.
Note that via the isomorphism
$\Hilb^n(\Bbb P((T_{(\Bbb P^2)^*})^*)/(\Bbb P^2)^*) \cong
\Hilb^n(F/(\Bbb P^2)^*)$, the fiber of the map $\pi$
over a point $[C] \in (\Bbb P^2)^*$ is isomorphic to $\Hilb^n(C)$.
So $\im(\alpha) = W_n$. Moreover, since $n \ge 2$,
an element $\xi \in \hil{n}$ is contained in at most one
line $C$ in $X=\Pee^2$. So if $\alpha(\xi') = \alpha(\eta')$
for $\xi', \eta' \in \Hilb^n(\Bbb P((T_{(\Bbb P^2)^*})^*)/(\Bbb P^2)^*)$,
then there exists a unique line $C$ such that
$\alpha(\xi') = \alpha(\eta') \subset C$, i.e.,
$\xi'$ and $\eta'$ are contained in the fiber $\pi^{-1}([C])$
over the point $[C] \in (\Bbb P^2)^*$. Since $\alpha$ maps
this fiber isomorphically to $\Hilb^n(C) \subset \hil{n}$,
we must have $\xi' = \eta'$.
By the Theorem 4.1.11 in \cite{Got},  
$\Hilb^n(\Bbb P((T_{(\Bbb P^2)^*})^*)/(\Bbb P^2)^*) \cong
\Bbb P(\Sym^n((T_{(\Bbb P^2)^*})^*))$. Thus, $\alpha$ gives the isomorphism
$$\Bbb P(\Sym^n((T_{(\Bbb P^2)^*})^*)) \cong W_n.  \tag 5.8$$

So the set of all the lines in
$\Hilb^n(C) \subset \hil{n}$ with $C$ being lines in $X=\Pee^2$ is
parameterized by the Grassmannian bundle
$\Grass(\Sym^n((T_{(\Bbb P^2)^*})^*), 2)$ over $(\Bbb P^2)^*$.
By Theorem 5.4,  from the universal properties of 
the Hilbert schemes and the Grassmannians,
we obtain a bijective morphism
$\Lambda: \Grass(\Sym^n((T_{(\Bbb P^2)^*})^*), 2) \to 
\frak M(\beta_{\ell}-(n-1)\beta_n)$.
Since both $\Grass(\Sym^n((T_{(\Bbb P^2)^*})^*), 2)$ and  
$\frak M(\beta_{\ell}-(n-1)\beta_n)$ are smooth,
the bijectivity of $\Lambda$ implies that $\Lambda$ is an isomorphism.
\endproof             

By Lemma 3.8 and (3.7), the morphism 
$\Phi_2=\frak P_2\circ \varphi_2\colon\,\, \hil{n} \to \Bbb P^{N_2}$
is an embedding associated to the complete linear 
system $|nD_{\ell}-B_n/2|$. By Theorem 3.14 (iii), 
$\gamma \subset \hil{n}$ is mapped to a line by $\Phi_2$ 
if and only if either $\gamma \sim \beta_n$
or $\gamma \sim (\beta_\ell - (n-1) \beta_n)$.
Therefore, regarding $\hil{n}$ as a closed subvariety of $\Bbb P^{N_2}$,
then the Hilbert scheme of lines in $\hil{n}$ is 
the disjoint union of $\frak M(\beta_n)$ and 
$\frak M(\beta_\ell - (n-1) \beta_n)$.

\section{6. The structure of the morphism $\varphi_1$}

Our goal in this final section is to study the structure of 
the morphism $\varphi_1$ defined in (3.5) when $n\ge 3$.
We recall from Proposition 3.12 (ii) that the structure of 
$\varphi_1$ has been determined when $n=2$. We begin with a lemma.

\lemma{6.1} Let $n \ge 3$, $X=\Pee^2$, and $W_n$ be 
the contraction locus of the morphism 
$\varphi_1: \hil{n} \to \varphi_1(\hil{n})$ defined by (3.5). Then, 

{\rm (i)} a fiber of $\varphi_1$ is either a point or
the subset $C^{(n)} \subset \hil{n}$ for some line $C \subset X$;

{\rm (ii)} the restriction $\varphi_1|_{W_n}: W_n \to \varphi_1(W_n)$
is isomorphic to the projective bundle
$\Bbb P(\Sym^n((T_{(\Bbb P^2)^*})^*))$ over $(\Bbb P^2)^*$,
where $T_{(\Bbb P^2)^*}$ is the tangent bundle of $(\Bbb P^2)^*$;

{\rm (iii)} $N \cong (\varphi_1|_{W_n})^*M \otimes \Cal O(-1)$ where
$N = N_{W_n \subset \hil{n}}$ is the normal bundle of $W_n$ in $\hil{n}$,
$M$ is some vector bundle over $\varphi_1(W_n) \cong (\Bbb P^2)^*$,
and $\Cal O(1)$ stands for the tautological line bundle over
the projective bundle $\varphi_1|_{W_n}: W_n \to \varphi_1(W_n)$.
\endstatement
\proof
(i) Fix $\xi \in \hil{n}$. By definition, 
$\varphi_1(\xi) \in \varphi_1(\hil{n}) \subset 
\Grass\big(H^0\big(X, \Cal O_X(n-1)\big), n\big)$ 
corresponds to the $n(n-1)/2$-dimensional subspace 
$H^0\big(X, I_{\xi}\otimes\Cal O_X(n-1)\big)$ in 
$H^0\big(X, \Cal O_X(n-1)\big)$. 
So $\eta \in (\varphi_1)^{-1}(\varphi_1(\xi))$ if and only if
$H^0\big(X, I_\eta \otimes\Cal O_X(n-1)\big) =
H^0\big(X, I_{\xi}\otimes\Cal O_X(n-1)\big).$
Let $\frak d \subset |(n-1)\ell|$ be the linear system corresponding to 
the subspace $H^0\big(X, I_{\xi}\otimes\Cal O_X(n-1)\big)$ in 
$H^0\big(X, \Cal O_X(n-1)\big)$. Then, either $\frak d$ cuts out 
a $0$-dimensional subscheme $\zeta$ with $\zeta \supset \xi$,
or $\frak d$ has a fixed component $C$.

Assume that $\frak d$ cuts out a $0$-dimensional subscheme $\zeta$ 
with $\zeta \supset \xi$. Then, 
$H^0\big(X, I_\zeta \otimes\Cal O_X(n-1)\big) =
H^0\big(X, I_{\xi}\otimes\Cal O_X(n-1)\big)$, and 
$(\varphi_1)^{-1}(\varphi_1(\xi)) \subset \{ \eta \in \hil{n}|\,\,
\eta \subset \zeta\}$. Conversely, 
if $\eta \in \hil{n}$ and $\eta \subset \zeta$,
then $H^0\big(X, I_\eta \otimes\Cal O_X(n-1)\big) \supset
H^0\big(X, I_{\zeta}\otimes\Cal O_X(n-1)\big) =
H^0\big(X, I_{\xi}\otimes\Cal O_X(n-1)\big)$.
So $H^0\big(X, I_\eta \otimes\Cal O_X(n-1)\big)=
H^0\big(X, I_{\xi}\otimes\Cal O_X(n-1)\big)$ 
since they have the same dimension. 
Thus $\eta \in (\varphi_1)^{-1}(\varphi_1(\xi))$,
and $\{ \eta \in \hil{n}|\,\, \eta \subset \zeta\} \subset 
(\varphi_1)^{-1}(\varphi_1(\xi))$. 
Hence $(\varphi_1)^{-1}(\varphi_1(\xi)) = \{ \eta \in \hil{n}|\,\,
\eta \subset \zeta\}.$
Now we claim $(\varphi_1)^{-1}(\varphi_1(\xi)) = \{ \xi \}$.
Indeed, if $(\varphi_1)^{-1}(\varphi_1(\xi))$ contains more than
one point, then there is an irreducible projective curve $\Gamma$ 
in $(\varphi_1)^{-1}(\varphi_1(\xi)) = 
\{ \eta \in \hil{n}|\,\, \eta \subset \zeta\}$ since the closed subset 
$\{ \eta \in \hil{n}|\,\, \eta \subset \zeta\}$ of $\hil{n}$ 
is known to be connected. Moreover, the projective curve $\Gamma$
is contracted to points by both the Hilbert-Chow morphism $\rho_n$
and the morphism $\varphi_1$. So $\Gamma \cdot D_\ell = 0$,
and $\Gamma \cdot B_n = 0$ in view of (3.7). By Theorem 3.14 (i), 
the intersection of $\Gamma$ with every divisor on $\hil{n}$ is zero. 
This is impossible. So $(\varphi_1)^{-1}(\varphi_1(\xi)) = \{ \xi \}$.

Next, assume that $\frak d$ has a fixed component $C$.
Let $C \sim d\ell$ where $d \ge 1$. Then, 
${n(n-1)/2} -1 = \dim (\frak d) \le \dim \, |(n-1-d)\ell| 
= {(n+1-d)(n-d)/2} -1.$
It follows that $d = 1$, $C$ is a line in $X = \Pee^2$, 
and $\xi \subset C$ (i.e. $\xi \in C^{(n)} \subset \hil{n}$).
So $(\varphi_1)^{-1}(\varphi_1(\xi)) \subset C^{(n)}$.
On the other hand, letting $p(C)$ be a linear homogeneous polynomial
defining the line $C$, we see that 
for every $\eta \in C^{(n)} \subset \hil{n}$,
the subspace $H^0\big(X, I_\eta \otimes\Cal O_X(n-1)\big)$ in 
$H^0\big(X, \Cal O_X(n-1)\big)$ consists of 
all the degree-$(n-1)$ homogeneous polynomials divisible by $p(C)$.
In particular, $H^0\big(X, I_\eta \otimes\Cal O_X(n-1)\big)=
H^0\big(X, I_\xi \otimes\Cal O_X(n-1)\big)$.
So $C^{(n)} \subset  (\varphi_1)^{-1}(\varphi_1(\xi))$.
Hence $(\varphi_1)^{-1}(\varphi_1(\xi)) = C^{(n)}$.

(ii) By (i), the contraction locus $W_n$ is the union of all the subsets
$C^{(n)} \subset \hil{n}$ with $C$ being lines in $X=\Pee^2$.
So (ii) follows from (5.8) and (i).

(iii) Fix a line $C$ in $X =\Pee^2$. By (ii), 
$C^{(n)}$ is a fiber in the fiber bundle
$\varphi_1|_{W_n}: W_n \to \varphi_1(W_n)$.
So we have $N_{C^{(n)} \subset W_n} = \Cal O_{C^{(n)}}^{\oplus 2}$.
Thus by Proposition 5.6 (i), we obtain an exact sequence of normal bundles:
$$\exact{\Cal O_{C^{(n)}}^{\oplus 2}}
{\Cal O_{C^{(n)}}^{\oplus 2}\oplus \Cal O_{C^{(n)}}(-1)^{\oplus (n-2)}}
{(N_{W_n \subset \hil{n}})|_{C^{(n)}} = N|_{C^{(n)}}}.$$
It follows that $N|_{C^{(n)}} \cong \Cal O_{C^{(n)}}(-1)^{\oplus (n-2)}$
for every line $C$ in $X = \Pee^2$.
Since these $C^{(n)}$ with $C$ being lines in $X = \Pee^2$
are precisely the fibers of the projective bundle
$\varphi_1|_{W_n}: W_n \to \varphi_1(W_n)$,
we conclude that $N \cong (\varphi_1|_{W_n})^*M \otimes \Cal O(-1)$
for some rank-$(n-2)$ vector bundle $M$ over
$\varphi_1(W_n) \cong (\Bbb P^2)^*$.
\endproof

Let $n \ge 3$. By Lemma 6.1 (ii), we have
$\dim(W_n) = (n+2) \le \dim(\hil{n})-1$.
So $\varphi_1: \hil{n} \to \varphi_1(\hil{n})$ is a birational morphism.
Since $K_{\hil{n}} = -3 D_\ell$ where $\ell$ is a line in $X=\Pee^2$, 
$(\beta_\ell - (n-1) \beta_n) \cdot K_{\hil{n}}= -3 < 0$.
This says that $R \,\, {\overset \text{def} \to =} \,\,
\Bbb R^+ \cdot (\beta_\ell - (n-1) \beta_n)$ is an extremal ray
since $(-K_{\hil{n}})$ is nef and the Picard number of $\hil{n}$ is $2$. 
Thus $\varphi_1$ is the contraction of $R$ since 
$((n-1)D_\ell-B_n/2) \cdot (\beta_\ell - (n-1) \beta_n) =0$. 
Moreover, when $n \ge 4$, $\dim(W_n) \le \dim(\hil{n})-2$, 
so $\varphi_1$ is a small contraction. 

\theorem{6.2} Let $n \ge 4$ and $X=\Pee^2$. Then, the flip for 
the small contraction $\varphi_1: \hil{n} \to \varphi_1(\hil{n})$ 
defined via (3.5) exists.
\endstatement
\proof
We shall prove that the Flip Conjecture 5-1-10 in \cite{KMM}
holds for $\varphi_1$. More precisely, we construct
the commutative diagrams (6.8) and (6.9) satisfying
this Flip Conjecture.
The proof is a little long, but may be divided into two steps.
In step one, we construct $\W \varphi_-$ and $\W \varphi_+$. 
In step two, we construct $\varphi_+$.

Let $\W \varphi_-: \W {\hil{n}} \to \hil{n}$ be the blowing-up
of $\hil{n}$ along the contraction locus $W_n$ of $\varphi_1$,
and $\W {W_n} \subset \W {\hil{n}}$ be the exceptional divisor.
Let $N$ be the normal bundle of $W_n$ in $\hil{n}$.
Then over $\W {W_n} =\Pee(N^*)$, we have the tautological surjection
$$(\W \varphi_-|_{\W {W_n}})^*(N^*) \to \Cal O_{\W {W_n}}(1) \to 0
\tag 6.3$$
where $\Cal O_{\W {W_n}}(1)$ is the tautological line bundle over 
$\W {W_n} =\Pee(N^*)$. By Lemma 6.1 (iii), 
$N \cong (\varphi_1|_{W_n})^*M \otimes \Cal O_{W_n}(-1)$
where $M$ is some rank-$(n-2)$ vector bundle over
$\varphi_1(W_n) \cong (\Bbb P^2)^*$,
and $\Cal O_{W_n}(1)$ stands for the tautological line bundle over
the projective bundle $\varphi_1|_{W_n}: W_n \to \varphi_1(W_n)$.
So from (6.3), we obtain a surjection:
$$(\varphi_1|_{W_n} \circ \W \varphi_-|_{\W {W_n}})^*(M^*) \to 
(\W \varphi_-|_{\W {W_n}})^*\Cal O_{W_n}(-1) \otimes 
\Cal O_{\W {W_n}}(1) \to 0.$$
By the Proposition 7.12 in Chapter II of \cite{Har}, the map
$(\varphi_1|_{W_n} \circ \W \varphi_-|_{\W {W_n}}):
{\W {W_n}} \to \varphi_1(W_n)$
factors through the natural projection $\Pee(M^*) \to \varphi_1(W_n)$. 
In fact, the induced morphism $\W \varphi: {\W {W_n}} \to \Pee(M^*)$ is 
a $\Pee^n$-bundle. So the fibers of 
$(\varphi_1|_{W_n} \circ \W \varphi_-|_{\W {W_n}}):
{\W {W_n}} \to \varphi_1(W_n)$
are naturally $\Pee^n \times \Pee^{n-3}$. Moreover, 
for the projective spaces $\Pee^n$ in the fibers $\Pee^n \times \Pee^{n-3}$,
we have $\Cal O_{\W {\hil{n}}}(\W {W_n})|_{\Pee^n}
\cong \Cal O_{\Pee^n}(-1)$, using the facts that
$\Cal O_{\W {\hil{n}}}(\W {W_n})|_{\Pee^n \times \Pee^{n-3}} 
\cong \Cal O(a, -1)$ for some integer $a$ and that
$$\align
 \Cal O_{\Pee^n}^{\oplus (n-2)} \otimes \Cal O_{\Pee^n}(1)
 \cong N^*|_{\Pee^n}
&\cong (\W \varphi_-|_{\Pee^n \times \Pee^{n-3}})_* \left [
 \Cal O_{\W {\hil{n}}}(-\W {W_n})|_{\Pee^n \times \Pee^{n-3}} \right]\\
&\cong \Cal O_{\Pee^n}^{\oplus (n-2)} \otimes \Cal O_{\Pee^n}(-a).\\
\endalign$$
In other words, the restriction 
$\Cal O_{\W {\hil{n}}}(-\W {W_n})|_{\W {W_n}}$ is a tautological line 
bundle of $\W \varphi: {\W {W_n}} \to \Pee(M^*)$.
By the results in \cite{Nak, F-N}, there exists a morphism
$\W \varphi_+: \W {\hil{n}} \to (\hil{n})_+$ such that 
$(\hil{n})_+$ is smooth and contains $\Pee(M^*)$.
Moreover, $\W \varphi_+|_{\W {W_n}} = \W \varphi$ and 
$\W \varphi_+$ is the blowing-up of $(\hil{n})_+$ along $\Pee(M^*)$.

Now we claim that $(\W \varphi_-)^*((n-1)D_\ell - B_n/2) =
(\W \varphi_+)^*H_+$ for some divisor $H_+$ on $({\hil{n}})_+$.
Indeed, recall from (3.7) that  
$((n-1)D_\ell - B_n/2) = (\varphi_1)^*H  $
for some very ample divisor $H$ on $\varphi_1(\hil{n})$.
Consider the divisor 
$$(\W \varphi_-)^*((n-1)D_\ell - B_n/2)
= (\varphi_1 \circ \W \varphi_-)^*H \tag 6.4$$
on $\W {\hil{n}}$. Since $\Zee \cdot [\W {W_n}] \oplus 
(\W \varphi_-)^* \text{Cl}({\hil{n}}) = \text{Cl}(\W {\hil{n}}) =  
\Zee \cdot [\W {W_n}] \oplus (\W \varphi_+)^* \text{Cl}(({\hil{n}})_+)$
where $\text{Cl}( \cdot )$ stands for the divisor class group,
we conclude that 
$$(\W \varphi_-)^*((n-1)D_\ell - B_n/2) =
\lambda \W {W_n} + (\W \varphi_+)^*H_+ \tag 6.5$$
for some integer $\lambda$ and some divisor $H_+$ on $({\hil{n}})_+$. 
Let $\W C \subset \W {\hil{n}}$ be a line contained in
a fiber $\Pee^n$ of the projective bundle
$\W \varphi: {\W {W_n}} \to \Pee(M^*)$.
Then, $\W \varphi_+(\W C)=\W \varphi (\W C)$ is a point in
$\Pee(M^*) \subset ({\hil{n}})_+$,
and $(\varphi_1 \circ \W \varphi_-)(\W C)$ is a point in  
$\varphi_1(W_n) \subset \varphi_1(\hil{n})$.
So applying $\cdot \W C$ to both sides of (6.5) and using (6.4),
we obtain
$$\align
0 &= (\varphi_1 \circ \W \varphi_-)^*H \cdot \W C
   = (\W \varphi_-)^*((n-1)D_\ell - B_n/2) \cdot \W C    \\
  &= (\lambda \W {W_n} + (\W \varphi_+)^*H_+) \cdot \W C
   = \lambda \cdot (\W {W_n} \cdot \W C)               
   = \lambda \cdot (c_1(\Cal O_{\Pee^n}(-1)) \cdot \W C) = -\lambda
\endalign$$
where we have used the fact that $\Cal O_{\W {\hil{n}}}(\W {W_n})|_{\Pee^n}
= \Cal O_{\Pee^n}(-1)$ when $\Pee^n$ is a fiber of
the projective bundle $\W \varphi: {\W {W_n}} \to \Pee(M^*)$.
So $\lambda = 0$ and 
$$(\W \varphi_-)^*((n-1)D_\ell - B_n/2) =(\W \varphi_+)^*H_+.
\tag 6.6$$

Since $\W \varphi_+$ is the blowing-up morphism of the smooth variety
$({\hil{n}})_+$, we see that $|(\W \varphi_+)^*H_+|
= (\W \varphi_+)^*|H_+|$. So combining with (6.4) and (6.6), we obtain
$$|(\varphi_1 \circ \W \varphi_-)^*H|
= |(\W \varphi_-)^*((n-1)D_\ell - B_n/2)|    
= |(\W \varphi_+)^*H_+| = (\W \varphi_+)^*|H_+|.  \tag 6.7$$
Since $\varphi_1$ is the contraction of the extremal ray
$R = \Bbb R^+ \cdot (\beta_\ell - (n-1) \beta_n)$, 
it is known from \cite{KMM} that 
$\varphi_1(\hil{n})$ has at most rational singularities.
So $|(\varphi_1 \circ \W \varphi_-)^*H|
=(\varphi_1 \circ \W \varphi_-)^*|H|$. In view of (6.7), 
$(\varphi_1 \circ \W \varphi_-)^*|H|
=|(\varphi_1 \circ \W \varphi_-)^*H|= (\W \varphi_+)^*|H_+|$.
Since $H$ is a very ample divisor on $\varphi_1(\hil{n})$,
$(\varphi_1 \circ \W \varphi_-)^*|H|$ is base-point-free.
It follows from $(\W \varphi_+)^*|H_+| =
(\varphi_1 \circ \W \varphi_-)^*|H|$ that 
$|H_+|$ is also base-point-free. So $|H_+|$ induces a morphism
$\varphi_+: ({\hil{n}})_+ \to \varphi_+(({\hil{n}})_+)$.
Moreover, since $(\W \varphi_+)^*|H_+| =
(\varphi_1 \circ \W \varphi_-)^*|H|$, 
we conclude that $\varphi_1 \circ \W \varphi_-
= \varphi_+ \circ \W \varphi_+$ up to an isomorphism between 
$\varphi_1(\hil{n})$ and $\varphi_+(({\hil{n}})_+)$.
So putting $\varphi_- = \varphi_1$, we obtain a commutative diagram:
$$\matrix
       &                    &\W {\hil{n}}       &        &           \\
    &\swarrow^{\W \varphi_-}&    &\searrow^{\W \varphi_+}&           \\
\hil{n}&                    &                   &       &(\hil{n})_+.\\
       &\searrow^{\varphi_-}&       &\swarrow^{\varphi_+}&           \\
       &                    &\varphi_1(\hil{n}) &        &           \\
\endmatrix  \tag 6.8$$ 
Note that $\varphi_+: (\hil{n})_+ - \Pee(M^*) \to 
\varphi_1(\hil{n}) - \varphi_1(W_n)$ is an isomorphism. 

It is standard to check that $K_{(\hil{n})_+}$ is $\varphi_+$-ample.
So the commutative triangle
$$\matrix
\hil{n}&&{\overset {\W \varphi_+ \circ (\W \varphi_-)^{-1}} 
        \to \dasharrow}&       &(\hil{n})_+                    \\
       &\searrow^{\varphi_-}&       &\swarrow^{\varphi_+}&     \\
       &                    &\varphi_1(\hil{n}) &        &     \\
\endmatrix  \tag 6.9$$ 
satisfies the Flip Conjecture 5-1-10 in \cite{KMM}, i.e.,
the flip for $\varphi_- = \varphi_1$ exists.
\endproof
\noindent
{\bf Remark 6.10.} The proof of Theorem 6.2
shows that for $n=3$, $\varphi_1: \hil{n} \to \varphi_1(\hil{n})$
is the blowing-up along $\varphi_1(W_n) \cong 
(\Pee^2)^*$, and $W_n$ is the exceptional locus.


\end{document}